\documentclass[11pt,article,reqno]{amsart}
\usepackage{amsmath, amsthm, amsfonts}

\usepackage[english]{babel}
\usepackage[square,numbers]{natbib}
\usepackage[normalem]{ulem} 

\usepackage{hyperref}

\usepackage{amsmath, amsthm, amsfonts}

\usepackage{amssymb,url}
\usepackage{mathtools}
\usepackage{amsaddr}
\usepackage[top=29mm,bottom=29mm,left=30mm,right=30mm]{geometry}
\usepackage{mathrsfs}
\usepackage[final]{graphicx}
\usepackage{ulem}
\usepackage{color}

\usepackage{mathtools}
\usepackage{graphicx}
\makeatletter
\def\section{\@startsection{section}{1}%
\z@{1\linespacing\@plus\linespacing}{1\linespacing}%
{\bf\centering}}
\def\subsection{\@startsection{subsection}{0}%
\z@{\linespacing\@plus\linespacing}{\linespacing}%
{\bf}}
\def\subsubsection{\@startsection{subsubsection}{0}%
\z@{\linespacing\@plus\linespacing}{\linespacing}%
{\bf}}
\makeatother

\makeatletter
\@addtoreset{equation}{section}
\makeatother

\newtheorem{theorem}{Theorem}[section]
\newtheorem{corollary}[theorem]{Corollary}
\newtheorem{lemma}[theorem]{Lemma}
\newtheorem{proposition}[theorem]{Proposition}

\newtheorem{remark}[theorem]{Remark}
\theoremstyle{definition}
\newtheorem{definition}[theorem]{Definition}
\newtheorem{example}[theorem]{Example}
\newtheorem{assumption}[theorem]{Assumption}

\def\al{\alpha}
\def\be{\beta}

\begin{document}
\title{scaling limit of dependent random walk}
\author{ Jeonghwa Lee}

\address{Department of Mathematics and Statistics, University of North Carolina Wilmington, USA}
\email{leejb@uncw.edu}

\thanks{
 \\ \medskip
\noindent
2020 {\it MS Classification}:  60F17, 60G07, 35D30}

\begin{abstract}
Recently, a generalized Bernoulli process (GBP) was developed as a stationary binary sequence that can have long-range dependence. In this paper, we find the scaling limit of a random walk that follows GBP.  The result is a new class of non-Markovian diffusion processes.  The limiting processes include continuous-time stochastic processes with stationary increments whose correlation decays with an exponential rate, a power law, or an exponentially tempered power law. The limit densities solve a time tempered fractional diffusion equation or time fractional diffusion equation. The second-family of Mittag-Leffler distribution and exponential distribution arise as special cases of the limiting distributions. Subordinated processes are considered as time-changed L\'evy processes, and the governing equations and dependence structure of the subordinated processes are discussed.  
\end{abstract}

\maketitle

\baselineskip 0.5 cm

\bigskip \medskip

\section{Introduction}

 A random walk describes the path of a walker who takes a random step at discrete, equally distanced time points. The scaling limit of a simple random walk is Brownian motion which was first formulated by \cite{einstein_1905_ber}   to describe the movements of particles in a liquid. 
 More generally, a L\'evy process is infinitely divisible with independent stationary increments, and it is considered a limit of a random walk that takes i.i.d. random steps, which may have an infinite second moment, at equal time intervals \citep{ jklafter_2011_first}. If 
a random walk can take large steps with a power-law jump distribution, the resulting L\'evy process can model super-diffusion where particles spread faster than the normal diffusion, 
and its limit density solves a space fractional diffusion equation. 

In a continuous time random walk (CTRW),  random steps are taken at random time points, and it can model both supper- and sub-diffusion \citep{meerschaert_2019_stochastic}.
 When the wait times between steps follow a power law distribution, the resulting limiting process is non-Markovian, governed by a (space) time fractional diffusion equation, and it 
 has a wide range of applications in finance, biology, hydrology,  etc   \citep{baeumer_2007_fractional,benson_2000_application,  deng_2006_parameter,sca,schumer_2003_fractal, weiss_2007_anomalous}.  Random
walks where jumps or wait times follow an exponentially tempered power law distribution can model tempered fractional diffusion \citep{ ALRAWASHDEH2017892, meerschaert_2008_tempered, SABZIKAR201514}.  
This exponential tempering has both mathematical and practical advantages, and it is proved to be useful in geophysics and finance \citep{ carr_2002_the, SABZIKAR201514}. The tempered fractional Brownian motion is the tempered fractional derivative or integral of a
Brownian motion and has stationary increments that can exhibit semi-long range dependence \citep{MEERSCHAERT20132269}.   

 More recently, a random walk that takes correlated steps was developed in \cite{leonenko_2018_correlated, meerschaert_2009_correlated}. 
In \cite{meerschaert_2009_correlated}, random steps follow a stationary linear process, which results in the time-changed fractional Brownian motion and linear fractional stable motion as the limiting processes. In \cite{leonenko_2018_correlated}, random steps are obtained by Markov chains through the Bernoulli urn-scheme model and Wright–Fisher model, which results in fractional Pearson diffusion. With a different approach, non-Markovian diffusion processes were developed in \cite{mura_2008_nonmarkovian} by introducing a memory kernel in the forward drift equation.  

In this paper, we consider a correlated random walk where steps follow a generalized Bernoulli process (GBP).
The GBP was first developed in \cite{lee_2021_generalized} as a stationary binary sequence whose covariance function decays with a power law. The sum in a finite sequence of GBP defines the fractional binomial distribution, and its variance is asymptotically proportional to the length of the sequence to a fractional power. In \cite{Lee19052025}, it was found that the GBP and the fractional Poisson process have the same large-scale behavior. More specifically, the fractional binomial distribution and the fractional Poisson distribution, when scaled, have the same limiting distribution, which suggests that GBP is a discrete analog of the fractional Poisson process.
 As the Poisson process can approximate the Bernoulli process under certain conditions, the relationship between the GBP and the fractional Poisson process can be considered its counterpart when there is a long-range dependence. In \cite{lee_2024_on}, GBP was broadened to include various covariance functions.

In this work, we find the scaling limit of a dependent random walk that follows GBP and the dependence structure of the limiting process. The result is a new class of continuous-time stochastic processes with stationary increments whose correlation decays with exponentially tempered power-law. 
It turns out that the probability density functions of the limiting processes solve 
time tempered fractional diffusion equations.  We also find the limiting process of a random walk that follows a non-stationary GBP, and the limiting distribution includes the second-type Mittag-Leffler distribution and exponentially decaying distribution as special cases. Subordinated processes are considered with a L\'evy process whose time is changed by the limiting process of the dependent random walk, which leads to non-Markovian processes that can serve as versatile models for diffusion processes. The marginal density of the resulting processes solves various space-time tempered fractional diffusion equations. Unlike other approaches where the probability density of a diffusion equation was found using an inverse stable subordinator \cite{ALRAWASHDEH2017892} or forward drift with power kernel  \cite{mura_2008_nonmarkovian}, our method allows us to characterize the full stochastic structure of the processes with its dependence structure, not just the marginal density functions.  


The paper is organized as follows. In Section 2, we introduce GBP and a dependent random walk that follows GBP. In Section 3, we find the limiting processes of the dependent random walk, providing their moment generating function (mgf) and dependence structure. 
Subordinated processes are considered as time-changed L\'evy processes in Section 4, followed by the conclusion in Section 5. All the proofs are in Section 6.

Let us end the introduction with a few words about the notation. We will let
$i,i_1,i_2,i_3\cdots \in\mathbb{N},$ and $i_1<i_2<i_3<\cdots,$ unless stated otherwise. For $x\in \mathbb{R},$
$\lfloor x \rfloor$ denotes the largest integer smaller than or equal to $x,$ and $[n]=\{ 1,2,\cdots, n\}$ for $ n\in \mathbb{N}.$ For any two quantities, $a,b,$ that depend on another variable $x$, $a\sim b$ indicates asymptotic equivalence, i.e., $a/b\to 1$ as $x$ goes to $0$ or $\infty.$
  For a set $C\subset \mathbb{R},$ $|C|$ indicates the number of elements in $C$, and the ordered elements are denoted by  $C(1)=\min C, C(2)=\min \{  C\setminus\{C(1)\} \},$ and in general, the $k$-th smallest element in $C$ is denoted by $C(k)$. An indicator function is denoted by $I( \cdot)$, and generic constants are denoted by $c'$ and $ c''$ that can change from line to line. For notational convenience, we denote $P((\cap_i A_i)\bigcap ( \cap_j B_j))$ by $P(\cap_i A_i \cap_j B_j).$

\section{Dependent Random Walk}

We first develop a generalized Bernoulli process (GBP) similar to \cite{lee_2021_generalized,lee_2024_on}. A dependent random walk is then defined as 
a random walk whose steps are dependent and follow GBP.
Let $\{f_n\}$ be a sequence of functions defined on $\mathbb{N},$  and $b\in (0,1].$ 
\begin{definition} 
Define for a set $ B=\{i_0, i_1, \cdots,i_k\}\subset \mathbb{N} $ with $i_0< i_1< \cdots<i_k, $
\[  L_n(B)= \prod_{j=1,\cdots,k}f_n(i_j-i_{j-1}) .\]
If $B=\emptyset$, define $ L_n(B):=1/(b f_n(\lfloor nT\rfloor))$, and if $|B|=1, L_n(B):=1. $ 
\end{definition}

 \begin{definition}
 Define for disjoint sets, $B,F\subset \mathbb{N}$ with $ |F|>0,$  \begin{align*}
    D_n(B,F)=
    \sum_{j=0}^{|F|}\sum_{ \substack{F'\subset F\\ |F'|=j}}(-1)^{j} L_n(B\cup F').
 \end{align*} 
    If $F=\emptyset, D_n(B,F):=L_n(B).$
    
    \end{definition}

  We define GBP in such a way that the dependence structure among ``1"s is induced by  $L_n(\cdot)$ using the distance between consecutive indices, and the dependence structure among ``1"s and ``0"s is derived by $D_n(\cdot, \cdot) $ using the inclusion-exclusion principle. 

  More specifically, for any fixed $T>0, n\in\mathbb{N},$ we define a sequence of dependent Bernoulli random variables $\{ \xi_{i}^{(n)}; i\in  [\lfloor nT \rfloor  ] \}$ with joint probabilities as follows. For any  $ B\subset [\lfloor nT \rfloor],$
  \begin{align}\label{p1}     
P(  \cap_{i\in B} (\xi_{i}^{(n)} =1 ) ) = b f_n( \lfloor nT \rfloor)  L_n(B),
  \end{align}
  i.e., for
  $1 \leq i_1<i_2<\cdots<i_k \leq \lfloor nT \rfloor,$ we have
\begin{align*}
 P(  \cap_{j\in [k]} (\xi_{i_j}^{(n)} =1 ) )=b f_n( \lfloor nT \rfloor) \prod_{j\in [k-1]} f_n(|i_{j+1}-i_j|) ,\end{align*} and the marginal probability becomes
$   P(  \xi_{i}^{(n)} =1)=b f_n( \lfloor nT \rfloor), i\in [\lfloor nT \rfloor].$
In general, by the inclusion and exclusion principle, the joint probabilities are defined as   \begin{align} \label{p2}
    P( \cap_{i'\in F }\{\xi_{i'}^{(n)}=0\}  \cap_{i\in B}   \{\xi_{i}^{(n)}=1\} ) &= \sum_{k=0}^{|F|}\sum_{\substack{F'\subset F\\ |F'|=k}} (-1)^{k}P(\cap_{i\in B\cup F'}\{\xi_{i}^{(n)}=1\} ) \nonumber\\&=b f_n( \lfloor nT \rfloor)D_n(B, F)
\end{align}
for any disjoint sets $B,F \subset [\lfloor nT \rfloor].$ In fact, the joint probabilities (2.1-2.2) are well-defined if and only if $D_n(B,F)\geq 0$ for any disjoint sets $B,F \subset [\lfloor nT \rfloor],$ which requires the following assumption on the function $f_n$.
\begin{assumption} \label{assumption}
{\it i)} $f_n: \mathbb{N} \to (0,1)$ is a decreasing function such that $f_n(x+1)/f_n(x) $ is non-decreasing as $x$ increases, and
\\ {\it ii)} 
    $f_n(2) >f_n(1)^2.$ 
\end{assumption}
 Assumption \ref{assumption} {\it i)} is satisfied for a decreasing function that decays as fast as or slower than exponential decay.
In the following examples, we show a class of functions that satisfy Assumption \ref{assumption}.  
\begin{example} 
{\it i)}$f_n(x)= c\exp{(-\lambda x/n)}$ where $c\in (0,1), \lambda >0.$ 
\\  {\it ii)} $f_n(x)= c x^{\alpha-1} $ for $c\in (0,2^{\alpha-1}), \alpha<1.$ 
\\  {\it iii)} 
\begin{align}\label{f}
    f_n(x)= c \exp{(-\lambda x/n)} x^{\alpha-1}  \end{align}
 where $ c\in (0,2^{\alpha-1}), \alpha \leq 1, \lambda >(\alpha-1)/T. $
\end{example}
The function $f_n$ determines the dependence structure of GBP through (\ref{p1}) and (\ref{p2}). In (\ref{f}), the parameters $\alpha$ and $\lambda$ affect the degree to which the correlation function of GBP decays with a power-law and an exponential rate, respectively. The parameter $c$ has a scaling effect on the limiting process of GBP. 

\begin{example}
{\it i)} $f_n(x)=p+c\exp(-x),$ where $p\in (0,1),$ and
\[   c\in \min\bigg\{ 1-p ,\frac{ (e^{-2}-2pe^{-1} ) + \sqrt{ ( e^{-2}-2pe^{-1}  )^2-4e^{-2}p(p-1)  } }{2e^{-2}} \bigg\}. \]
\\ {\it ii)} $f_n(x)=p+cx^{\alpha-1},$ where $ p\in (0,1),  \alpha <1,$ and
\[ c\in \min\bigg\{ 1-p ,0.5\Big({(2^{\alpha-1}-2p ) + \sqrt{ (2^{\alpha-1}-2p )^2-4p(p-1)  } }  \Big) \bigg\}   . \] 
\end{example}

The next theorem proves that under Assumption \ref{assumption}, $D_n(B,F)>0,$ for any disjoint sets $B,F \subset [\lfloor nT \rfloor],$ therefore, the sequence of dependent Bernoulli variables is well defined with the joint probabilities (2.1-2.2).
\begin{theorem} \label{prop0} Under Assumption \ref{assumption}, a stationary binary sequence $\{\xi_i^{(n)}, i\in  [\lfloor nT \rfloor]\}$ is well-defined with probabilities (\ref{p1},\ref{p2}) and 
\[ \mathrm{corr} (\xi_i^{(n)},  \xi_{i+k}^{(n)}  )=\frac{ f_n(k)-b f_n ( \lfloor nT \rfloor) }{ 1-b f_n ( \lfloor nT \rfloor) } ,\] for $k\neq0,$ and $ i,i+k\in  [\lfloor nT \rfloor].$
\end{theorem}
We call this stationary process a generalized Bernoulli process (GBP) and denote it by  \[ \{\xi_i^{(n)}, i\in  \lfloor nT \rfloor\} \sim GBP(f_n).\] In the rest of the paper, we specifically focus on GBP$(f_n)$ with $f_n$ in (\ref{f}),  derive its limiting processes, and find the dependence structures and governing equations of the limiting processes. Different functions $f_n$, such as the one in Example 2.5 $\it ii)$,  require different approaches than what is used in this paper to find the limiting process, which will be future work. 
\begin{assumption} \label{real_assum}
    $f_n$ is of the form (\ref{f}) with $c\in (0,2^{\alpha-1}), \alpha \in (0,1], \lambda >(\alpha-1) /T.$ 
\end{assumption}

\begin{remark}

 Let $   \xi_i^{(n)},i=1,2,\cdots, \lfloor nT \rfloor,$ be GBP($f_n$)  satisfying Assumption \ref{real_assum}. Then
\[ \sum_{k=1}^{n} \mathrm{corr}( \xi_1^{(n)},  \xi_{1+k}^{(n)}  ) \sim\begin{dcases}
c_1 n   &\text{ if }b>0, \alpha=1, \\    c_2  n^{\alpha}
 &\text{ if }  \alpha\in (0,1)  ,  \end{dcases}    \] 
where  $c_1,c_2$ are some constants. The results can be interpreted that the $GBP(f_n)$ has a strong correlation that decays slowly due to either a long-range dependence from the power function in (\ref{f}) or shrinking distance between variables (i.e., sampling more densely in a fixed domain) by the scaling factor $n$ in (\ref{f}) or both. 
\end{remark}

\section{Limiting Theorems}

In this section, we will show the scaling limit of a dependent random walk that follows GBP. More specifically, we will show the weak convergence of  a scaled random walk, \begin{equation}\label{sc1}
X_t^{(n)}= \frac{\sum_{i\in [\lfloor nt \rfloor]  } \xi_i^{(n)}  }{ n^{\alpha} } \hspace{8pt} \text{   for } 0<t \leq T, \end{equation}  in the space $D[0,T]$ with the Skorokhod ($J_1$)  topology, where $  \xi_i^{(n)} $ follows $ $GBP$(f_n)$ satisfying Assumption \ref{real_assum}.

\begin{theorem} \label{theo1} Let $\{\xi_i^{(n)}, i\in [\lfloor nT \rfloor]\}$ be GBP($f_n$) with $f_n$ satisfying Assumption \ref{real_assum}. 
Then, 
\[   X_t^{(n)}\Rightarrow X_t\]  in the space $D[0,T]$ with the Skorokhod ($J_1$)  topology. The limiting process $X_t$ has stationary increments, is continuous in probability, and the moment generating function (mgf) of its marginal distribution is  
\begin{align} \label{mgf1}
    \mathrm{E}(\exp(sX_t))=m(s,t)= 1+ 
  t a_T s + a_T  \sum_{k=1}^{\infty} s^{k+1} \nu^{k}    
 \int_{0}^t \int_{0}^x   \frac{e^{-\lambda y}y^{k\alpha-1}}{\Gamma(k\alpha)} dy    dx  , \end{align} where $a_T=cb e^{-\lambda T}T^{\alpha-1}$ and $\nu=c\Gamma(\alpha)$.
The mgf of any finite-dimensional distribution of $X_t$  is as follows. 
For any $m\in \mathbb{N}, 0=t_0<t_1<t_2<\cdots <t_m \leq T,$ and $\{s_j\in \mathbb{R}; j\in [m]\},$
\begin{align}\label{mgf2}
& \mathrm{E}(\exp(s_1X_{t_1}+s_2X_{t_2}+\cdots +s_mX_{t_m} )) 
:=m((S_1,S_2,\cdots,S_m),(t_1,t_2,\cdots,t_m))\nonumber \\&=1 +a_T \sum_{\substack{C\subset [m] \\ C\neq \emptyset  }}\sum_{\substack{k_i\in\mathbb{N}\\ i\in  [|C|] }}  t_{C(|C|)}^{\alpha(\sum_{i\in  [|C|] }k_i-1)+1 }  \Big(\prod_{i\in {[|C|]}}S_{C(i)}^{k_i}\Big)A({C,\{k_i\} })
\end{align}
where  $S_i=\sum_{j=i}^ms_j,$  and
\begin{align*}
  A( C, \{k_i\} )=\iint \cdots \int_{Q_C}
 & \Big\{\frac{\nu^{k_1-1}(x_1-x_1')^{\alpha(k_1-1)-1}}{\Gamma((k_1-1)\alpha ) } \exp(-\lambda t_{C(|C|)} (x_1-x_1'))\Big\}^{I(k_1>1)} \\&\times \prod_{i\in C_2,i\geq2} 
 \Big\{ c ( x_{i}'-x_{i-1})^{\alpha-1} \frac{ \nu^{k_i-1} (x_{i}-x_{i}')^{\al(k_i-1)-1}}{\Gamma((k_i-1)\al)}\Big\}\\ &\times
 \prod_{j\in C_1} 
 \Big\{ c ( x_{j}-x_{j-1})^{\alpha-1} \Big\} \exp\big({-\lambda t_{C(|C|)} (x_{|C|}-x_1)}\big)
  \prod_{i\in C_2}dx_i' \prod_{\ell\in [|C|]}dx_{\ell}  ,
\end{align*}
where
\begin{align*}
 Q_C=\Big\{ x_i',x_i, x_j:  
\frac{ t_{C(i)-1}}{t_{C(|C|)}}
 <x_i'<x_i<\frac{t_{C(i)}}{t_{C(|C|)}}, i\in C_2,  
\frac{ t_{C(j)-1}}{t_{C(|C|)}}
 <x_j<\frac{t_{C(j)}}{t_{C(|C|)}}, j\in C_1
 \Big\}, \end{align*}
 $C_1=\{ i: 1\leq i\leq |C|, k_i= 1\}$ and  
$C_2=\{ i: 1\leq i\leq |C|, k_i\geq 2\}. $ 
\end{theorem}

From (\ref{mgf1}, \ref{mgf2}), we can find the moments of $X_t$ and its dependence structure as follows. For any $m\in \mathbb{N}, 0=t_0<t_1<t_2<\cdots <t_m \leq T,$ $C\subset[m], \{k_i\in \mathbb{N}; i\in [|C|]\},$ we have
\[ \mathrm{E}\Big(\prod_{i\in [|C|]} (X_{t_{C(i)}}-X_{t_{C(i)}-1})^{k_i} \Big)=a_Tt_{C(|C|)}^{\alpha(\sum_{i\in  [|C|] }k_i-1)+1 } A(C,\{k_i\})\prod_{i\in[|C|]}k_i!. \]
For example, if $C=\{2,5,6\},$ and $k_1=k_2=k_3=1,$
\begin{align*}
   \mathrm{E}[( X_{t_2}-X_{t_{1}}) ( X_{t_5}-X_{t_{4}})( X_{t_6}-X_{t_{5}})]=a_Tt_6^{2\al+1}\int_{t_5/t_6}^1\int_{t_4/t_6}^{t_5/t_6}\int_{t_1/t_6}^{t_2/t_6}c(x_3-x_2)^{\al-1}c(x_2-x_1)^{\al-1}\\
   \times \exp({-\lambda t_6(x_3-x_1)}) dx_1dx_2dx_3 .
\end{align*}
For $C=\{1,3,6\},$ and $k_1=3, k_2=4,k_3=2,$
\begin{align*}
   &\mathrm{E}[X_{t_1}^3 ( X_{t_3}-X_{t_{2}})^4( X_{t_6}-X_{t_{5}})^2]=3!4!2!a_Tt_6^{8\al+1}\int_{t_5/t_6}^1\int_{t_5/t_6}^{x_3}\int_{t_2/t_6}^{t_3/t_6}\int_{t_2/t_6}^{x_2}\int_{0}^{t_1/t_6}\int_{0}^{x_1}\frac{\nu^{2}(x_1-x_1')^{2\al-1}}{\Gamma(2\al)} \\ &\times
   (x_2'-x_1)^{\al-1}\frac{\nu^{3}(x_2-x_2')^{3\al-1}}{\Gamma(3\al)} (x_3'-x_2)^{\al-1}\frac{\nu(x_3-x_3')^{\al-1}}{\Gamma(\al)}
    \exp({-\lambda t_6(x_3-x_1')}) dx_1'dx_1dx_2'dx_2dx_3'dx_3 .
\end{align*}
It is found that  the correlation of increments has an exponentially tempered power law decay, i.e., it decays with a power law initially but eventually decays exponentially fast for a large lag.
\begin{corollary} \label{cor1}
{\it i)} $ \mathrm{E}(X_t)=t a_T, $ and for $k\geq 2,$
\[ \mathrm{E}(X_t^k)= t^{\alpha(k-1) +1}\frac{ a_T k! \nu^{k-1}}{ \Gamma((k-1)\alpha)} \int_0^1 \int_0^{x} e^{-\lambda ty} y^{\alpha(k-1)-1}   dydx .\] 
{\it ii)} For $0<t<s\leq T,$
\begin{align*}
    \mathrm{E}(X_t X_s)=& s^{\alpha+1} c{a_T} \int_{t/s}^1\int_{0}^{t/s} (x_2-x_1)^{\al-1} e^{-\lambda s(x_2-x_1) }  dx_1 dx_{2} \\&+ t^{\alpha+1} {2ca_T} \int_0^1 \int_{0}^x e^{-\lambda ty} y^{\alpha-1}  dy dx, \text{    and}   
\end{align*}

 \begin{align*}
    \mathrm{cov}(X_t, X_s)=& s^{\alpha+1} c{a_T} \int_{t/s}^1\int_{0}^{t/s} (x_2-x_1)^{\al-1} e^{-\lambda s(x_2-x_1) }  dx_1 dx_{2} \\&+ t^{\alpha+1} {2ca_T} \int_0^1 \int_{0}^x e^{-\lambda ty} y^{\alpha-1}  dy dx\\& -  ts    a_T^2 .
\end{align*} 
{\it iii)}
Let $ \nabla X_{t}= X_{t}-X_{t-d}$ be the increment of $X_t.$  For any $s,t>d>0,$ 
and as $d=\Delta\to 0,$
\begin{align*}
\mathrm{E}(\nabla X_s \nabla X_{s+t})&\sim  c'\Delta^2  e^{-\lambda t} t^{\alpha-1} , \\ \mathrm{cov}(\nabla X_s, \nabla X_{s+t})&\sim c' \Delta^2 (e^{-\lambda t} t^{\alpha-1}- b e^{-\lambda T}T^{\alpha-1}) , \end{align*} 
where $ c'= a_Tc$,
and 
\[ \mathrm{corr}(\nabla X_s, \nabla X_{s+t})\sim \frac{e^{-\lambda t} t^{\alpha-1}- b e^{-\lambda T}T^{\alpha-1} }{  
    c'' \Delta^{\alpha-1}- b e^{-\lambda T}T^{\alpha-1} },\] where $c''=2/(\alpha(\alpha+1)).$ 
\\ {\it iv)}
  For any $m\in \mathbb{N},$ $0<t_1<t_2<\cdots<t_m\leq T, $ and $k_i \in \mathbb{N},i=1,2,\cdots,m,$ 
    \[ \lim_{d\to 0}\frac{\mathrm{E}\big(\prod_{i=1}^m (\nabla X_{t_i})^{k_i}\big)}{\prod_{i=1}^m \mathrm{E}\big((\nabla X_{t_i})^{k_i}\big) 
 }    =  e^{-\lambda (t_m-t_1)}(c/a_T)^{m-1}\prod_{i=2}^m  (t_i-t_{i-1})^{\alpha-1} . \]
    \end{corollary}

Next, we define a different version of scaled random walk, 
 \begin{equation}\label{sc2}
    X_t^{(*,n)} =\frac{\sum_{i\in [\lfloor nt \rfloor]  } \xi_i^{(*,n)}  }{ n^{\alpha} }\hspace{8pt} \text{ for } t>0,\end{equation}
where $ \{\xi_i^{(*,n)} , n \in \mathbb{N}  \} $ is the nonstationary GBP defined as follows.

For any set $B \subset \mathbb{N}$, let   \begin{align} \label{p3}
    &P(  \cap_{i\in B}   \{\xi_i^{(*,n)} =1\} ) = P(  \cap_{i\in B_{+1}}   \{\xi_i^{(n)} =1\} )/P(\xi_{1} =1) ,
\end{align} where $B_{+1}:= \{1,i+1: i\in B \} $ (if $B=\emptyset, B_{+}:=\{1\})$ and $\{ \xi_i^{(n)}, i\in [\lfloor nT\rfloor] \} \sim GBP (f_n)$ for  $T> (\max B +1)/n.$ Note that the probability (\ref{p3}) is independent of $T$ and  $b$. 
For any disjoint sets $B,F \subset \mathbb{N},$ we will let
\begin{align} \label{p4}
    &P( \cap_{i'\in F }\{\xi_{i'}^{(*,n)}=0\}  \cap_{i\in B}   \{\xi_i^{(*,n)} =1\} ) = \sum_{k=0}^{|F|}\sum_{\substack{ F' \subset F\\ |F'|=k}} (-1)^{k}P(\cap_{i\in F'\cup B}\{\xi_i^{(*,n)} =1\} ),
\end{align}
  with which we will show that $\{ \xi_{i}^{(*,n)}, i\in \mathbb{N} \}$ is  well defined for any $n\in \mathbb{N}$, and it can be considered as a binary sequence in GBP which starts after the first ``1" appears.

\begin{proposition} \label{prop_01} Under Assumption \ref{real_assum}, a binary sequence $\{\xi_i^{(*,n)}, i\in  \mathbb{N}\} $ is well-defined with probabilities (\ref{p3},\ref{p4}). It is non-stationary with  $  P(\xi_i^{(*,n)}=1)=f_n(i)$ for any $i\in \mathbb{N},$ and 
\[ \mathrm{cov} (\xi_i^{(*,n)},  \xi_{j}^{(*,n)}  )={ f_n(i)f_n(j-i)- f_n (i)f_n (j) }\] for  $ i<j.$ We call it GBP*$(f_n)$.
\end{proposition}

\begin{theorem} \label{theo2} Let  $ \{ \xi_i^{(*,n)} , i\in \mathbb{N}\} $ be  GBP*$(f_n)$ with $f_n$ satisfying Assumption  \ref{real_assum}. Then, 
\[   X_t^{(*,n)}\Rightarrow X_t^*\]  in the space $D[0,\infty)$ with the Skorokhod ($J_1$)  topology. The limiting process $X_t^*$ is continuous in probability, and the mgf of its marginal distribution is    \begin{align} \label{mgf3}
\mathrm{E}(\exp(sX_t^*))=m^*(s,t)=1+\sum_{k=1}^{\infty}s^k \nu^{k} \int_{0}^t \frac{e^{-\lambda y}y^{k\alpha-1}}{\Gamma(k\alpha)} dy  , \end{align} where $\nu=c\Gamma(\alpha)$.  The mgf of  a finite-dimensional distribution is, for any $m\in \mathbb{N}, 0=t_0<t_1<t_2<\cdots <t_m,$ and $ s_i\in \mathbb{R}, i\in [m]$, 
\begin{align} \label{mgf4}
&\lim_{n\to \infty} \mathrm{E}(\exp(s_1X_{\alpha,t_1}^*+s_2X_{\alpha,t_2}^*+\cdots +s_m X_{\alpha,t_m}^*)):=m^*((S_1,S_2,\cdots,S_m),(t_1,t_2,\cdots,t_m)) \nonumber \\&=1 +\sum_{\substack{C\subset [m] \\ C\neq \emptyset  }}\sum_{\substack{k_i\in\mathbb{N}\\ i\in [ |C| ]}} t_{C(|C|)}^{\alpha\sum_{i\in [|C|]}k_i } \Big(\prod_{i\in [|C| ] }S_{C(i)}^{k_i}\Big) A^*( C, \{k_i\} )  ,
\end{align}
where 
$S_i=\sum_{j=i}^ms_j,$ and

\begin{align*}
  A^*( C, \{k_i\} )=\iint \cdots \int_{Q_C^*}
   \prod_{i\in C_2} 
 \Big\{ c ( x_{i}'-x_{i-1})^{\alpha-1} \frac{ \nu^{k_i-1} (x_{i}-x_{i}')^{\al(k_i-1)-1}}{\Gamma((k_i-1)\al)}\Big\}\\ \times
 \prod_{j\in C_1} 
 \Big\{ c ( x_{j}-x_{j-1})^{\alpha-1} \Big\} \exp\big({-\lambda t_{C(|C|)} x_{|C|}}\big)
  \prod_{i\in C_2}dx_i' \prod_{\ell\in [|C|]}dx_{\ell}  ,
\end{align*}
where
\begin{align*}
 Q_C^*=\Big\{ x_i',x_i, x_j:  x_0=0,
\frac{ t_{C(i)-1}}{t_{C(|C|)}}
 <x_i'<x_i<\frac{t_{C(i)}}{t_{C(|C|)}}, i\in C_2,  
\frac{ t_{C(j)-1}}{t_{C(|C|)}}
 <x_j<\frac{t_{C(j)}}{t_{C(|C|)}}, j\in C_1
 \Big\}, \end{align*}
 $C_1=\{ i: 1\leq i\leq |C|, k_i= 1\}$ and  
$C_2=\{ i: 1\leq i\leq |C|, k_i\geq 2\}. $ 
\end{theorem}

From (\ref{mgf3}, \ref{mgf4}), for any $m\in \mathbb{N}, 0=t_0<t_1<t_2<\cdots <t_m \leq T,$ $C\subset[m], \{k_i\in \mathbb{N}; i\in [|C|]\},$ we have
\[ \mathrm{E}\Big(\prod_{i\in [|C|]} (X_{t_{C(i)}}^*-X_{t_{C(i)}-1}^*)^{k_i} \Big)=t_{C(|C|)}^{\alpha\sum_{i\in  [|C|] }k_i } A^*(C,\{k_i\})\prod_{i\in[|C|]}k_i!. \]
For example, if $C=\{2,5,6\},$ and $k_1=k_2=k_3=1,$
\begin{align*}
   \mathrm{E}[( X_{t_2}^*-X_{t_{1}}^*) ( X_{t_5}^*-X_{t_{4}}^*)( X_{t_6}^*-X_{t_{5}}^*)]=t_6^{3\al}\int_{t_5/t_6}^1\int_{t_4/t_6}^{t_5/t_6}\int_{t_1/t_6}^{t_2/t_6}c(x_3-x_2)^{\al-1}c(x_2-x_1)^{\al-1}cx_1^{\al-1}\\
   \times \exp({-\lambda t_6x_3}) dx_1dx_2dx_3 .
\end{align*}
For $C=\{1,3,6\},$ and $k_1=3, k_2=4,k_3=2,$
\begin{align*}
   &\mathrm{E}[(X_{t_1}^*)^3 ( X_{t_3}^*-X_{t_{2}}^*)^4( X_{t_6}^*-X_{t_{5}}^*)^2]=3!4!2!t_6^{9\al}\int_{t_5/t_6}^1\int_{t_5/t_6}^{x_3}\int_{t_2/t_6}^{t_3/t_6}\int_{t_2/t_6}^{x_2}\int_{0}^{t_1/t_6}\int_{0}^{x_1}(x_1')^{\al-1}\frac{\nu^{2}(x_1-x_1')^{2\al-1}}{\Gamma(2\al)} \\ &\times
   (x_2'-x_1)^{\al-1}\frac{\nu^{3}(x_2-x_2')^{3\al-1}}{\Gamma(3\al)} (x_3'-x_2)^{\al-1}\frac{\nu(x_3-x_3')^{\al-1}}{\Gamma(\al)}
    \exp({-\lambda t_6x_3}) dx_1'dx_1dx_2'dx_2dx_3'dx_3 .
\end{align*}

\begin{corollary} \label{cor2}
   {\it i)} For $k\geq 1,$
 \[ \mathrm{E}\big( (X_t^*)^k \big)=t^{k\alpha}\frac{k!\nu^k}{\Gamma(k\alpha)}\int_0^1 e^{-\lambda ty} y^{k\alpha-1}dy  .\]
For $0<t<s,$
\begin{align*}
    \mathrm{E}(X_s^* X_t^*)&= s^{2\alpha } { c^2}\int_{t/s}^{1}\int_{0}^{t/s } x_1^{\al-1}(x_2-x_1)^{\alpha-1} \exp({-\lambda s x_{2}})  dx_1 dx_{2} \\&+ t^{ 2\alpha} \frac{2\nu^2}{\Gamma(2\alpha)} {\int_0^1 e^{-\lambda ty} y^{2\alpha-1} dy},
\end{align*} and
\begin{align*}
   \mathrm{cov}(X_s^*, X_t^*)&= s^{2\alpha } { c^2}\int_{t/s}^{1}\int_{0}^{t/s } x_1^{\al-1}(x_2-x_1)^{\alpha-1} \exp({-\lambda s x_{2}})  dx_1 dx_{2} \\&+ t^{ 2\alpha} \frac{2\nu^2}{\Gamma(2\alpha)} {\int_0^1 e^{-\lambda ty} y^{2\alpha-1} dy}\\&-(ts)^{\alpha} c^2\int_0^1 e^{-\lambda ty} y^{\alpha-1} dy \int_0^1 e^{-\lambda sx} x^{\alpha-1} dx.
\end{align*}

{\it ii)}
Let $ \nabla X_{t}^*= X_{t}^*-X_{t-d}^*$ be the increment of $X_t^*.$  For any $s,t>d>0,$  as $d=\Delta\to 0,$
\begin{align*}
\mathrm{E}(\nabla X_s^* \nabla X_{s+t}^*)&\sim  c'\Delta^2  e^{-\lambda (s+t)} s^{\al-1}t^{\alpha-1} , \end{align*} 
where $ c'= c^2$.
For any $m\in \mathbb{N},$ $0=t_0<t_1<t_2<\cdots<t_m, $ and $k_i \in \mathbb{N},i=1,2,\cdots,m,$  
\begin{align*}
  \mathrm{E}\Big(\prod_{i=1}^m (\nabla X_{t_i}^*)^{k_i}\Big) \sim & \Delta^{\sum_{i\in[m]}(\alpha (k_i-1)+1)} \exp(-\lambda t_m)  \prod_{i=1}^m \Big\{ 
  \frac{k_i! \nu^{k_i-1} }{\Gamma((k_i-1)\alpha+2 ) }  c(t_j-t_{j-1} )^{\al-1}   \Big\}, 
\end{align*} and
    \[ \lim_{d\to 0}\frac{\mathrm{E}\big(\prod_{i=1}^m (\nabla X_{t_i}^*)^{k_i}\big)}{\prod_{i=1}^m \mathrm{E}\big((\nabla X_{t_i})^{k_i}\big) 
 }    =  e^{-\lambda t_m} (c/a_T)^{m}\prod_{i=1}^m  (t_i-t_{i-1})^{\alpha-1} . \]
\end{corollary} 
From the mgfs in (\ref{mgf1}, \ref{mgf2}, \ref{mgf3}, \ref{mgf4}), it is observed that the constant $c$ has a scaling effect on $X_t, X_t^*.$ Without loss of generality, we will let $c=1/\Gamma(\alpha),$ which results in 
   $\nu=1$ in the mgfs, call the corresponding stochastic processes   
the exponential Mittag-Leffler (EML) processes for reasons that will be found in the following sections, and denote them by $X_t\sim$ EML$( \alpha, \lambda ,b,T)$ and $X_t^*\sim$ EML$^*( \alpha, \lambda ).$
We also use the notations $X_t( \alpha, \lambda , b, T)$ and $X_t^*( \alpha, \lambda   )$ for the EML processes, and 
$ p(x,t)$ and $ p^*(x,t) $ for their marginal probability density functions (pdfs), respectively.
It turns out that the marginal distributions of the EML processes solve  tempered time-fractional diffusion equations. 

Let $D_t^{\alpha,\lambda }$  be the Riemann–Liouville tempered fractional derivatives defined in \cite{  can} and \cite{ SABZIKAR201514}.
   \[ D_t^{\alpha,\lambda } f(t) =\frac{e^{-\lambda t}}{\Gamma(n-\alpha)} \frac{d^n}{dt^n}\int_{0}^t (t-s)^{-\alpha+n-1} { e^{\lambda s}f(s)} ds , \hspace{8pt} n-1\leq \alpha < n,\] 
   for $n\in \mathbb{N}.$ If $ \lambda=0,$  $ D_t^{\alpha,0 } $ becomes the Riemann-Liouville  fractional
derivative, which we also denote by $D^{\al}_t$. If  $ \alpha=1,$  \[D_t^{1,\lambda } f(t)= e^{-\lambda t} \frac{ d ( e^{\lambda t} f(t)) }{dt}.\]

\begin{proposition} \label{prop}
For the EML processes with parameters, $\lambda \geq 0, \alpha\in(0,1],$ the following holds for their marginal pdfs. \\ 
{\it i)} The pdfs  solve tempered fractional diffusion equations,
    \begin{align*}
      D_t^{\alpha,\lambda }p^*(x,t) &=-\frac{\partial p^*(x,t)}{\partial x} +\delta(x)  D_t^{\alpha,\lambda }(1), \\  D_t^{\alpha,\lambda }p(x,t)& =-\frac{\partial p(x,t)}{\partial x}+ \Big(\delta(x)   D_t^{\alpha,\lambda }(1)+ \delta'(x) -\delta'(x){a_T}   D_t^{\alpha,\lambda }( t)   \Big), \end{align*}
    where
  $ \partial/ \partial x$ is a weak derivative, $\delta(x)$ is the Dirac delta function, and $\delta'(x)$ is the first distribution derivative of the Dirac delta function.
  \\ {\it ii)} The pdf of  $X_t^*(\alpha, \lambda )$ is
   \[p^*(x,t)= e^{-\lambda t}h(x,t)+\lambda \int_0^t e^{-\lambda y} h(x,y) dy\]
for $t,x\geq0,$ where $h$ is the pdf of inverse stable subordinator, more specifically, \[ h(x,t)=\frac{t}{x\alpha} g(t,x),\] where $g$ is the pdf of stable subordinator whose Laplace transform is  \[ \int_0^{\infty} e^{-st} g(t,x) dt = e^{-xs^{\alpha}}.\] The pdf of $X_t(\alpha,\lambda)$ is
    \[ p(x,t)= - a_T \int_0^t \frac{\partial p^*(x,y)}{\partial x} dy +\delta(x),\] for $x\geq 0, t\in [0,T].$ 
\end{proposition}

\begin{remark}
{\it i)} The EML$^*( \alpha, \lambda )$ process $X_t^*$ is defined for $t \in\mathbb{R}_+,$
whereas EML$( \alpha, \lambda ,b,T)$  process $X_t$ is defined for  $t\in [0,T],$ for any $T>0.$ One can use different parametrization for $X_t$, and denote
$\{ X_t; 0<t\leq T\} \sim $EML($ \alpha, \lambda , a_T $)
where $0<a_T\leq e^{-\lambda T}T^{\alpha}/\Gamma(\alpha).$ The latter can avoid the parameter identification problem 
since  all finite-dimensional distributions are uniquely defined with these parameters. With the first parametrization,
$X_t\sim$ EML$( \alpha, \lambda ,b_1,T_1)$
and $Y_t\sim $EML$( \alpha, \lambda ,b_2,T_2)$ have the same finite-dimensional distributions on $0<t<\min\{T_1, T_2\}$ 
if  $b_1 e^{-\lambda T_1}T_1^{\alpha}=b_2 e^{-\lambda T_2}T_2^{\alpha}. $
In this paper, however, we will use the original parametrization since it gives an insight on the role of $b$ and $T$ on the dependence structure of the EML process.
\\ {\it ii)} From (\ref{mgf2}) and (\ref{mgf4}), one can find (generalized) self-similar property in the EML processes, i.e., for any $ u>0,$

\[ \{  X_{ut} ( \alpha, \lambda , b, T) ; t\in [0,T/u] \} \stackrel{f.d.d}{=} \{ u^{\alpha}X_t ( \alpha, u\lambda , b, T/u) ; t\in [0,T/u]\}, \] and
\[ \{  X_{ut}^* ( \alpha, \lambda ) ; t>0 \} \stackrel{f.d.d}{=} \{ u^{\alpha}X_t^* ( \alpha, u\lambda ) ; t>0 \} .\]  
\\ {\it iii)} In the EML$^*(\al, \lambda)$ process, {when $\lambda>0,$}  
   the moments can be expressed with the incomplete gamma function, and (\ref{mgf3}) can be written as
   \[    \mathrm{E}(e^{sX_t^*})= 1+\sum_{k=1}^{\infty}  (\lambda ^{-\alpha}s)^k \frac{\gamma(k\alpha, \lambda t)}{\Gamma(k\alpha)}  , \]
   where $\gamma(k,t)= \int_0^t  e^{-x}x^{k-1} dx.$ When $ t\to \infty, $ the exponential-Mittag-Leffler variable $X^*_t$  converges in distribution to the exponential random variable with scale parameter $\lambda ^{-\alpha}$, since the mgf of $X_t^*$ converges to the mgf of the exponential distribution.
  \\ {\it iv)}
    In \cite{ALRAWASHDEH2017892}, it was shown that the pdf of an inverse
tempered stable subordinator solves a time tempered fractional diffusion equation. We remark that the pdf and the tempered fractional diffusion equation it solves in \cite{ALRAWASHDEH2017892} are not the same as those in Proposition \ref{prop} in our paper.  We take a different approach, using a dependent random walk and its scaling limit, to find a pdf that solves a tempered time-fractional diffusion equation, which allows us to find the complete stochastic structure of the diffusion process with its dependence structure. 
\end{remark}
\subsection{ Exponential Processes}
We will examine a special case in the EML process. More specifically, we will look at EML$^*(\alpha,\lambda )$ and  EML$(\alpha,\lambda ,b,
T)$ when $\alpha=1$ and $\lambda >0$. By Theorems \ref{theo1} and \ref{theo2},  for GBP with  $f_n(x)= c\exp{(-\lambda x/n)},$  $c\in (0,1), \lambda >0,$ the scaled random walk
$ X_t^{(n)} $ in (\ref{sc1}) converges weakly to  $cX_t(1,\lambda , b,T)$ in the space $D[0,T]$, and 
$X_t^{(*,n)}$ in (\ref{sc2}) converges weakly to $cX_t^*(1,\lambda )$ in the space $D[0,\infty],$ with the Skorokhod   topology.  The  marginal distribution of the limiting process $X_t(1,\lambda , b,T)$ 
has mgf
\begin{align} \label{mgf5} \mathrm{E}(\exp(sX_t))= m(s,t)&
=1+{t b e^{-\lambda T} s}+b e^{-\lambda T}s^2\Big(\frac{t-(\lambda -s)^{-1}(1-e^{-t(\lambda -s)})}{\lambda -s}\Big),   
\end{align}  
and the marginal distribution of $X_t^*(1,\lambda )$ has mgf
\begin{align} \label{mgf6}
 \mathrm{E}( \exp(sX_t^* ))=m^*(s,t)=1+s\Big(\frac{1-e^{-t(\lambda -s)}}{\lambda -s}\Big).\end{align}

From Proposition \ref{prop}, we have the following results on the pdfs of $X_t^*(1,\lambda )$ and  $X_t(1,\lambda , b,T).$  The pdf of $X_t^*$ is
\begin{align} \label{exp_eq}
 p^*(x,t)=I_{(0,t]}(x){\lambda }e^{ -\lambda x}+\delta_{t}(x)e^{-\lambda t} ,\end{align}
 where $\delta_a(x)=\delta(x-a)$ is the Dirac delta function, and
\[  \frac{\partial p^*(x,t)}{\partial t}+\lambda p(x,t) = -\frac{\partial p^*(x,t)}{\partial x}+\lambda \delta(x), \] from which one obtains 
the pdf of  $X_t$ as 
\begin{align*}
 p(x,t)& =-be^{-\lambda T}\int_0^t \frac{d p^*(x,y)}{dx} dy +\delta(x)  \\&=  
I_{(0,t]}(x) b   e^{-\lambda T} {\lambda (2+
\lambda t-\lambda x)}e^{ -\lambda x}+\delta_{t}(x)b   e^{-\lambda T}e^{-\lambda t}+\delta(x)
(1-b e^{-\lambda T}-t\lambda b   e^{-\lambda T}).
\end{align*}
The pdf of $X_t^*$ is exponentially decaying until $x=t$ with a point mass, and $X_t$ has point masses at 0 and $t.$
As $t \to \infty,$ $X_t^*$ converges in distribution to the exponential random variable.  
We call these limiting processes exponential processes, and denote them by  $\{X_t; t\in [0,T] \} \sim $ Exp$(\lambda,b,T)$ and $\{ X_t^*; t>0 \}\sim$  Exp$^*(\lambda)$.

From the mgfs in (\ref{mgf5}, \ref{mgf6}), we obtain
\begin{align*}
\mathrm{E}(  X_t )= tb  e^{-\lambda T}, \hspace{8pt}
\mathrm{E}(  X_t^{*} )=  \frac{1-e^{-\lambda t}}{\lambda },\end{align*}
and for $k\geq 2,$
\begin{align*}
    \mathrm{E}(X_t^{*k})&=k!  \int_{0}^t e^{-\lambda y} \frac{y^{k-1}}{(k-1)!} dy,\\
    \mathrm{E}(X_t^k)&
    =k! b e^{-\lambda T} \int_{0}^t \int_0^y e^{-\lambda x} \frac{x^{k-2}}{(k-2)!} dx dy .
\end{align*} 
 The Exp$(\lambda,b,T)$ process has stationary increments whose correlation decays exponentially fast. 
\begin{proposition} \label{prop2}
 Let $\nabla X_t=X_t-X_{t-d}$ be the increment of the exponential process, Exp$(\lambda,b,T)$. \\
{\it i)} For $0<t<s<T $ such that $ 0<d<s-t,$ 
\begin{align*}
    \mathrm{E}(\nabla X_t \nabla X_s)=c'e^{-\lambda (s-t)}(e^{\lambda d}-2+e^{-\lambda d}  ),
\end{align*}
where $c'= b e^{-\lambda T}/\lambda ^2.$
As $d=\Delta \to 0,$
\begin{align*}
    \mathrm{E}(\nabla X_t \nabla X_s)&\sim  \Delta^2 c' e^{-\lambda (s-t)},\\
   \mathrm{cov}(\nabla X_t, \nabla X_s)&\sim   \Delta^2  c' ( e^{-\lambda (s-t)}-b e^{-\lambda T}),
\end{align*}
where $c'=b e^{-\lambda T},$ and
\begin{align*}
   \mathrm{corr}(\nabla X_t, \nabla X_s)\sim \frac{   e^{-\lambda (s-t)}-b e^{-\lambda T} }{1 -b e^{-\lambda T}}.
\end{align*}
{\it ii)} 
If $0<s-t<d,$ as $s-t\to 0,$
\begin{align*}
   \mathrm{corr}(\nabla X_t, \nabla X_s)\sim 1- c' (s-t), 
\end{align*}
where $c'= 2b e^{-\lambda T} /(\lambda^2var(\nabla X_t)).$
\\{\it iii)} For any $m\in \mathbb{N},$ $0<t_1<t_2<\cdots<t_m\leq T, $ and $k_i \in \mathbb{N},i=1,2,\cdots,m,$ as $\Delta\to 0,$
\[  \mathrm{E}\Big(\prod_{i=1}^m (\nabla X_{t_i})^{k_i}\Big)    \sim \Delta^{\sum_{i\in [m]} k_i} c' e^{-\lambda (t_m-t_1)}, \]
and 
    \[ \lim_{d\to 0}\frac{\mathrm{E}\big(\prod_{i=1}^m (\nabla X_{t_i})^{k_i}\big)}{\prod_{i=1}^m \mathrm{E}\big((\nabla X_{t_i})^{k_i}\big) 
 }    = (1/c')^{m-1} e^{-\lambda (t_m-t_1)}, \] where $c'=be^{-\lambda T}. $
\end{proposition}

\begin{proposition}
 Let $\nabla X_t^*=X_t^*-X_{t-d}^*$ be the increments of  Exp$^*(\lambda)$. \\
{\it i)} For $0<t<s $ such that $ 0<d<s-t,$ 
as $\Delta\to 0,$
\begin{align*}
    \mathrm{E}(\nabla X_t^* \nabla X_s^*)&\sim  \Delta^2 e^{-\lambda s}.
\end{align*}
 For any $m\in \mathbb{N},$ $0<t_1<t_2<\cdots<t_m, $ and $k_i \in \mathbb{N},i=1,2,\cdots,m,$   
\[  \mathrm{E}\Big(\prod_{i=1}^m (\nabla X_{t_i}^*)^{k_i}\Big)    \sim  \Delta^{\sum_{i\in [m]} k_i}  e^{-\lambda t_m}. \]

\end{proposition}

\begin{remark}
We showed that the exponential process $\{cX_t^*(1,\lambda) ,t\geq 0\}$ is the limiting process of the dependent random walk that follows GBP $\{\xi_i^{(*,n)}, i\in \mathbb{N}\}$ with $f_n(x)=c\exp(-\lambda x/n)$ for $c\in (0,1).$ In fact, one can show similarly that GBP is well defined when $c=1$ (i.e. $f_n(x)=\exp(-\lambda x/n)),$  and the limiting process is $X_t^*(1,\lambda)$ with the pdf (\ref{exp_eq}). In the latter, the GBP can freeze, i.e., by (\ref{p2}) it can be shown that once ``0" appears, everything followed will be ``0" with probability 1, which explains the exponential decay in the pdf in (\ref{exp_eq}).  To explain this heuristically, note that for $i'< \lfloor nt \rfloor $
\[ f_n(i')=P(\text{time when first ``0" appears in } \{\xi_i^{(*,n)} \}  >i') ,\]
\[ f_n(i'-1)-f_n(i')=P( \text{time when first ``0" appears} =i' )=P(X_t^{(*,n)}=i'/n ) , \]
from which it follows that
\[ \frac{1}{n}\lambda \exp(-\lambda x)\sim P(X_t^{(*,n)}=x ) . \]
In the former case of $c<1,$ the GBP does not freeze, but  given $ \xi_{i'}^{(*,n)}=1,$  $ \{ \xi_i^{(*,n)} , i<i'\}$ become independent Bernoulli variables for any $i'\in \mathbb{N}$.   
\end{remark}

\subsection{ Mittag-Leffler Processes}
 We will investigate the EML$^*(\alpha, \lambda )$ and  EML$(\alpha, \lambda ,b, T)$ when $\lambda =0.$ The following results are readily derived from Theorems \ref{theo1} and \ref{theo2}. 

For GBP with  $f_n(x)= cx^{\alpha-1},$ where $c\in (0,1),\alpha\in (0,1),$ the scaled random walk
$ X_t^{(*,n)} (X_t^{(n)}) $ converges weakly to  $ \nu X_t^*(\alpha, 0)$ $   (\nu X_t(\alpha, 0, b,T))$
in the space $D[0,\infty)(D[0,T])$ with the Skorokhod  topology, where $\nu=c\Gamma(\alpha)$.    The marginal distribution of $X_t^*(\alpha, 0)$ is the second family of  Mittag-Leffler distribution, i.e.,  
\[\mathrm{E}(\exp(sX_t^*))=m(s,t)=E_{\alpha,1}(s t^{\alpha}), \] 
 and the mgf of
 $X_t(\alpha, 0,b, T)$ is  
\begin{align*}
\mathrm{E}( \exp(sX_t ))&:= m(s,t)=1+  s c' t E_{\alpha, 2}(s t^{\alpha} ), \end{align*} where $c'=b  T^{\alpha-1}/\Gamma(\alpha),$ and    \[E_{\alpha,\beta}(z)= \sum_{k=0}^{\infty}\frac{z^k}{\Gamma(k\alpha +\beta)} \] is the Mittag-Leffler function.
We will call  these limiting process  the Mittag-Leffler processes, and denote them by  $\{X_t; t\in[0,T]\}\sim$ ML$(\alpha, b, T)$  and  $\{X_t^*; t>0\}\sim$ ML$^*(\alpha)$. 
By Proposition \ref{prop}, it follows that their marginal pdfs solve time-fractional diffusion equations,
\begin{align} \label{mitag_pdf}
 \partial_t^{\alpha}{p}^*(x,t)&= - \frac{\partial}{\partial x} {p}^*(x,t)  ,\\
\partial_t^{\alpha}{p}(x,t)&=-  \frac{\partial}{\partial x} {p}(x,t)+ \delta'(x) \Big(1-\frac{c' t^{1-\alpha}}{\Gamma(2-\alpha)} \Big) \nonumber  , \end{align}  where $\partial_t^{\alpha}$ is the Caputo fractional derivative,
\[\partial_t^{\alpha}f(t)= \frac{1}{\Gamma(1-\alpha)} \int_0^t f'(x) (t-x)^{-\alpha} dx, \hspace{8pt} \alpha\in (0,1).    \]  
The solution to the equation (\ref{mitag_pdf}) is the pdf of the inverse stable subordinator \citep{mura_2008_nonmarkovian},
\[  p^*(x,t)= {t^{-\alpha} M_{\alpha} ( x t^{-\alpha} )},  \hspace{8pt} x,t\geq 0,\] where
\[  M_{\alpha}(x)=\sum_{k=0}^{\infty} \frac{(-x)^k}{k! \Gamma(-\alpha k+(1-\alpha) )},   \hspace{8pt} x\geq 0.\] 
Also, by Proposition \ref{prop} {\it ii}, \[p(x,t)= - c' \int_0^t \frac{\partial p^*(x,y)}{\partial x} dy +\delta(x), \] for $x\geq0, t\in[0,T], $ and $c'=bT^{\al-1}/\Gamma(\al).$   

Using the mgf of finite-dimensional distributions (\ref{mgf4}) with $\lambda=0,$ one can show that the ML$^*(\al)$ process is self-similar. For any $u>0,$
\[ \{  X_{ut}^*  ; t>0 \} \stackrel{f.d.d}{=} \{ u^{\alpha}X_t^*  ; t>0\} .\]
By Corollaries \ref{cor1} and \ref{cor2},  we can obtain the following properties of ML processes.  
For $k\in \mathbb{N},$ \begin{align*}
     \mathrm{E}(X_{t}^k )= c_{k,1} t^{\alpha (k-1)+1}, \hspace{8pt}
    \mathrm{E}((X_{t}^*)^k)= c_{k,2}  t^{\alpha k},\end{align*} 
        where $c_{k,1}= k!b T^{\alpha-1} /(\Gamma(\alpha)\Gamma(\alpha (k-1)+2)) $ and $c_{k,2}=k!/\Gamma(\alpha k+1) .$
   For $0<t<s \leq T, $
   \[   \mathrm{E}(X_{t}X_{s})=\frac{c' }{\Gamma(\alpha+2)}\big( {t}^{\alpha+1 }+{s}^{\alpha+1 }-(s-t)^{\al+1}\big)\]
   \[\mathrm{cov}(X_{t},X_{s} )= \frac{c' }{\Gamma(\alpha+2)}\big( {t}^{\alpha+1 }+{s}^{\alpha+1 }-(s-t)^{\al-1}\big)-(c')^2 ts, \]
    \[\mathrm{cov}(X_{t}^*,X_{s}^* )=\bigg( {t}^{2\alpha } \frac{2}{\Gamma(2\alpha+1)} +{s}^{2\alpha }\frac{1}{\Gamma(\al)^2}\int_{t/s}^1 \int_0^{t/s}  x_1^{\al-1}(x_2-x_1)^{\al-1} dx_1 dx_2  \bigg)-\frac{1}{ \Gamma(\alpha+1)^2} (st)^{\alpha},\]
where $ c'=b T^{\alpha-1}/\Gamma(\alpha).$ The ML$(\alpha, b, T)$ process has stationary increments whose correlation decays with a power law.

\begin{proposition} 
\label{prop3}
  Let  
$\nabla X_t=X_t-X_{t-d}$ be the increment of  the ML$(\alpha, b, T)$.
\\ {\it i)} For $ s-t>d,$
 \[\mathrm{E}( \nabla X_{ t} \nabla X_{s})= c' \big((s-t+d)^{\alpha+1}-2(s-t)^{\alpha+1}+ (s-t-d)^{\alpha+1} \big)  , \]
 where $c'=b T^{\alpha-1} /(\Gamma(\alpha)\Gamma(\alpha+2)).$ As $d= \Delta\to 0,$
 \begin{align*}
    \mathrm{E}(\nabla X_{t} \nabla X_{s})&\sim \Delta^2 c'(s-t)^{\alpha-1}   ,\\
\mathrm{cov}(\nabla X_{t}, \nabla X_{s})&\sim \Delta^2 c' \big((s-t)^{\alpha-1}-b T^{\alpha-1}\big)  ,\end{align*} where $ c'= {b T^{\alpha-1} /\Gamma(\alpha)^2 } , $
 and 
\[\mathrm{corr}(\nabla X_{t}, \nabla X_{s})\sim \frac{ (s-t)^{\alpha-1}-b T^{\alpha-1}}{ c''
\Delta^{\alpha-1}-b T^{\alpha-1}},  \]
where 
$c''= 2/(\alpha (\alpha+1)).$
\\ {\it ii)} For  $0<s-t<d,$ as $s-t \to 0,$
\[\mathrm{corr}(\nabla X_{t}, \nabla X_{s}) \sim 1- c'(s-t)^{\alpha+1} ,\]
where $c'=2bT^{\al-1}/(\Gamma(\al)\Gamma(\al+2)var(\nabla X_{t} ))
.$
\\ 
{\it iii)}
 For any $m\in \mathbb{N},$ $0<t_1<t_2<\cdots<t_m\leq T, $ and $k_i \in \mathbb{N},i=1,2,\cdots,m,$ 
    \[ \lim_{d\to 0}\frac{\mathrm{E}\big(\prod_{i=1}^m (\nabla X_{t_i})^{k_i}\big)}{\prod_{i=1}^m \mathrm{E}\big((\nabla X_{t_i})^{k_i}\big) 
 }    =  \prod_{i=2}^m  \frac{(t_i-t_{i-1})^{\alpha-1}}{c'}   ,\]
 where $c'=bT^{\al-1}.$
 \end{proposition}

\begin{proposition}
Let $ \nabla X_{t}^*= X_{t}^*-X_{t-d}^*$ be the increment of ML$^*(\al).$  For any $s,t>d>0,$  as $d=\Delta\to 0,$
\begin{align*}
\mathrm{E}(\nabla X_s^* \nabla X_{s+t}^*)&\sim  c'\Delta^2   s^{\al-1}t^{\alpha-1} , \end{align*} 
where $ c'= \Gamma(\al)^{-2}$.
For any $m\in \mathbb{N},$ $0=t_0<t_1<t_2<\cdots<t_m, $ and $k_i \in \mathbb{N},i=1,2,\cdots,m,$  
\begin{align*}
  \mathrm{E}\Big(\prod_{i=1}^m (\nabla X_{t_i}^*)^{k_i}\Big) \sim & \Delta^{\sum_{i\in[m]}(\alpha (k_i-1)+1)}   \prod_{i=1}^m \Big\{ 
  \frac{k_i!  }{\Gamma((k_i-1)\alpha+2 )\Gamma(\al) } (t_j-t_{j-1} )^{\al-1}   \Big\}. 
\end{align*} 
\end{proposition}

We note that the stochastic process that solves the time-fractional diffusion equation (\ref{mitag_pdf}) was found as an inverse stable subordinator in \cite{meerschaert_2019_stochastic}. Also, in \cite{mura_2008_nonmarkovian}, random processes that solves the equations (\ref{exp_eq},\ref{mitag_pdf}) were found by the forward drift equation with an exponential-decay kernel and power memory kernel. 
With our approach,  they were found as the scaling limits of dependent random walks, and the full stochastic structure of the processes was revealed with their dependence structure. 

\section{Subordinated Processes}

In this section, we examine the characteristic function and dependence structure of a subordinated process, a L\'evy process subordinated to EML processes. 
Let $Y(t)$ be L\'evy process with characteristic function $\mathrm{E}(e^{{\rm i}
sY(t)})=e^{t\psi(s)},$   $Y=Y(1),$ and
$ Z_t= Y(X_t),Z_t^*= Y(X_t^*)$ where $X_t$ and  $X_t^*$ are 
EML processes.

\begin{theorem} \label{theo3}
    The characteristic functions of finite-dimensional distributions of the subordinated processes, $Z_t, Z_t^*,$ are
 \[  \mathrm{E}(e^{{\rm i} s_1Z_{t_1}+ {\rm i}s_2 Z_{t_2}+\cdots + {\rm i}s_m Z_{t_m} })  =m\Big( \big( \psi (S_1) ,\psi (S_2) , \cdots ,\psi (S_m)  \big), 
(t_1,t_2,\cdots, t_m)   \Big),  \]
\[  \mathrm{E}(e^{{\rm i}s_1Z_{t_1}^*+ {\rm i}s_2 Z_{t_2}^*+\cdots +{\rm i} s_m Z_{t_m}^* })=m^*\Big( \big( \psi(S_1) ,\psi(S_2) , \cdots ,\psi(S_m)  \big), 
(t_1,t_2,\cdots, t_m)   \Big), \]
where $m, m^*$ are functions in (\ref{mgf2}) and (\ref{mgf4}), and $S_i=\sum_{j=i}^m s_j.$
\end{theorem}

Note that if  L\'evy process $Y$ has a finite mean, then by Corollaries 3.2 and 3.5,
\[ \mathrm{E}(Z_t^*)=\mathrm{E}(Y)\mathrm{E}(X_t^*)=\mathrm{E}(Y)t^{k\alpha}\frac{k!\nu^k}{\Gamma(k\alpha)}\int_0^1 e^{-\lambda ty} y^{k\alpha-1}dy ,\]
and 
\[ \mathrm{E}(Z_t)=\mathrm{E}(Y)\mathrm{E}(X_t)=\mathrm{E}(Y)ta_T.   \]

$Z_t$ is a stationary increment process, and if L\'evy process $Y$ has a finite mean,
\[ \mathrm{E}(  \nabla Z_t \nabla Z_s   ) =\mathrm{E}(Y)^2 \mathrm{E}( \nabla X_t  \nabla X_s ), \]
\[ \mathrm{cov}(  \nabla Z_t, \nabla Z_s   ) =\mathrm{E}(Y)^2 \mathrm{cov}( \nabla X_t,  \nabla X_s ) ,\]
where $\nabla X_t= X_t-X_{t-d}, \nabla Z_t= Z_t-Z_{t-d},$ and $t,s\in[d,T], s-t>d.$

\begin{corollary}\label{cor3} If L\'evy process $Y$ has a non-zero finite mean, $\mathrm{E}(Y)\neq 0,$ then
for the increment process $\nabla Z_t= Z_t-Z_{t-d}, $
 as $d=\Delta\to 0,$ 

{\it i)} 
\begin{align*}
\mathrm{E}(\nabla Z_s \nabla Z_{s+t})&\sim  c'\Delta^2  e^{-\lambda t} t^{\alpha-1} , \\ \mathrm{cov}(\nabla Z_s, \nabla Z_{s+t})&\sim c' \Delta^2 (e^{-\lambda t} t^{\alpha-1}- b e^{-\lambda T}T^{\alpha-1}) , \end{align*} 
for any $0<s,t\leq T,$ and
\\{\it iii)}    
    \[ \lim_{d\to 0}\frac{\mathrm{E}\big(\prod_{i=1}^m \nabla Z_{t_i}\big)}{\prod_{i=1}^m \mathrm{E}\big(\nabla Z_{t_i}\big) 
 }    =  e^{-\lambda (t_m-t_1)}(c/a_T)^{m-1}\prod_{i=2}^m  (t_i-t_{i-1})^{\alpha-1} , \]
   for any $m\in \mathbb{N},$ $0<t_1<t_2<\cdots<t_m\leq T. $
\end{corollary}

\begin{example}
 Let $Y(t)$ be stable L\'evy motion with index $ \beta\in (0,2],$ and 
\begin{align*} 
  \mathrm{E}(e^{{\rm i} sY})&=
  \exp\big[ -C' |s|^{\beta}(1+{\rm i}\gamma   \text{sgn}(s) \Phi  \big],
  \end{align*}
  where $ \gamma\in [-1,1], C'>0,$ and
  \[ \Phi=\begin{dcases} \tan \frac{\pi \beta}{2} & \text{ if } \beta\neq 1,\\ -\frac{2}{\pi}\log|s| & \text{ if } \beta=1 .\end{dcases}\] If $\beta\neq 1,$  we can write it as \begin{align} \label{psi}
  \mathrm{E}(e^{{\rm i} sY})&=\exp\big[ C ( p({\rm i} s )^{\beta} +q(-{\rm i} s )^{\beta}) \big]
  \end{align}
for a constant $C $ which is positive if $1<\beta\leq 2,$ negative if $0<\beta<1,$ and $p,q\geq0$ with $p+q=1.$  Let $X_t, X_t^*$ be ML processes, $X_t\sim$ML$(\al, b, T),$    $X_t^*\sim$  ML$^*(\alpha)$.  Then, the subordinated processes, $Z_t, Z_t^*,$ have the characteristic functions,
\begin{align*}
\mathrm{E}(e^{{\rm i}s Z_t })  &= 1+ (c't) C[ p({\rm i} s )^{\beta} +q(-{\rm i} s )^{\beta}]  E_{\alpha,2}(  C[ p({\rm i} s )^{\beta} +q(-{\rm i} s )^{\beta}]  t^{\alpha})  ,\\
 \mathrm{E}(e^{{\rm i}s Z_t^* }) & = E_{\alpha,1}(  C[ p({\rm i} s )^{\beta} +q(-{\rm i} s )^{\beta}]  t^{\alpha})  ,\end{align*}
where $c'=bT^{\al-1}/\Gamma(\al),$ and the marginal pdfs of $Z_t, Z_t^*$ solve  space-time fractional diffusion equations, 
\begin{align*}
   \partial_t^{\alpha}p(x,t)=&Cp D_x^{\beta}  p(x,t) +Cq D_{-x}^{\beta} p(x,t) -(1-a_TD_t^{\al}(t))(CpD_x^{\beta}\delta(x)+CqD_{-x}^{\beta}\delta(x)) ,  \\
  \partial_t^{\alpha}p^*(x,t)=&Cp D_x^{\beta}  p(x,t) +Cq D_{-x}^{\beta} p(x,t) ,   \end{align*} respectively.
The process $Z_t$ has stationary increments, and if $\beta\in (1,2]$ and $\mathrm{E}(Y)\neq 0,$   
\[ \mathrm{E}(  \nabla Z_s \nabla Z_{t}   ) \approx c'\Delta^2  (s-t)^{\alpha-1},\] for $0<t<s\leq T,$  and in general, for $0<t_1<t_2<\cdots <t_k\leq T,$
\[    \mathrm{E}(  \nabla Z_{t_1} \nabla Z_{t_2}\cdots \nabla Z_{t_k}   ) \approx c''  \Delta^k \prod_{i=2}^k (t_i-t_{i-1})^{\alpha-1}  ,\]
for some constants $c',c'',$ and $ \nabla Z_t= Z_t-Z_{t-d}$ for small $d=\Delta\approx 0$. 
\end{example}

\begin{example}
Let  $X_t,X_t^*$ be 
EML$(\alpha,\lambda,b,T)$, EML$^*(\alpha,\lambda ),$ respectively, and $Y(t)$ be tempered fractional  L\'evy motion with \[  \mathrm{E}(e^{{\rm i} s Y})=\exp\big[-C((\mu+{\rm i}s)^{\beta}-\mu^{\beta})\big]\] for $0<\beta<1$, $C>0$. The pdf of $Y(t)$ solves time
tempered fractional diffusion equation,
\[  \frac{\partial p(x,t)}{\partial t}=-C \partial_{x}^{\beta,\mu} p(x,t) , \]
where $\partial_{x}^{\beta,\mu}$ is the (positive) tempered fractional derivative of order $0<\beta<1,$
\[   \partial_{x}^{\beta,\mu} f(x)=  D_x^{\beta,\mu}f(x) -\mu^{\beta}f(x) ,\] defined in 
\cite{meerschaert_2019_stochastic}.

Then, the pdfs of $Z_t,Z_t^*$ solve space-time tempered fractional diffusion equations, 
\begin{align*}
   D_t^{\alpha,\lambda }p(x,t)&= C  \partial_{x}^{\beta,\mu} p(x,t)  +\Big(\delta(x)  D_t^{\alpha,\lambda }(1)-C  \partial_{x}^{\beta,\mu}\delta(x)\big(1-a_TD_t^{\al,\lambda}(t)\big) \Big), \\   D_t^{\alpha,\lambda }p^*(x,t)&=  C  \partial_{x}^{\beta,\mu} p^*(x,t)   +\delta(x)  D_t^{\alpha,\lambda }(1), \end{align*} respectively. 
 The subordinated process $Z_t$ has stationary increments
 $ \nabla Z_t= Z_t-Z_{t-d},$ and   for small $d=\Delta\approx 0,$
\begin{align*}
 \mathrm{E}(  \nabla Z_s \nabla Z_{t}   ) &\approx c'\Delta^2 e^{-\lambda (s-t)} (s-t)^{\alpha-1}, \hspace{8pt} 
\\
\mathrm{cov}(  \nabla Z_t, \nabla Z_s   ) &\approx c'\Delta^2 \Big( e^{-\lambda (s-t)} (s-t)^{\alpha-1} -b e^{-\lambda T} T^{\alpha-1} \Big) ,\end{align*}
for some constant $c'$ and $0<t<s\leq T$, and
\[    \mathrm{E}(  \nabla Z_{t_1} \nabla Z_{t_2}\cdots \nabla Z_{t_k}   ) \approx c''  \Delta^k e^{-\lambda(t_k-t_1)} \prod_{i=2}^k (t_i-t_{i-1})^{\alpha-1}  ,\]
for some constant $c''$ and $0<t_1<t_2<\cdots<t_k\leq T $. 
\end{example}
\section{Conclusion}
We found the scaling limits of random walks that take correlated steps and the dependence structure of the limiting processes, which leads to a new class of non-Markovian processes. The limiting processes include continuous-time stochastic processes with stationary increments whose correlation decays with an exponential rate, a power law, or an exponentially tempered power law.
We showed that the marginal densities of the limiting processes solve time tempered fractional diffusion equations or time fractional diffusion equations. When L\'evy processes are subordinated to the limiting processes, it results in various non-Markovian diffusion processes that are governed by space-time tempered fractional diffusion equations and can have stationary increments whose correlation decays with tempered power law.






\section{Proofs}

The following Lemmas 6.1 and 6.2 are similar to Lemmas 8.1 and 8.2 in \cite{lee_2024_on}.

\begin{lemma} \label{lem1}
Let $f_n: \mathbb{N} \to (0,1)$ be a decreasing function such that $f_n(x+1)/f_n(x)$ is non-decreasing as $x$ increases, then \[\frac{f_n(x+a)}{f_n(x) }\] is non-decreasing as $x$ increases for any fixed $a.$
\begin{proof}
We need to show for any $x,a\in \mathbb{N},$
\[ \frac{f_n(x+a)}{f_n(x) }\leq  \frac{f_n(x+1+a)}{f_n(x+1) }, \]
or 
\[ \frac{f_n(x+1)}{f_n(x) }\leq \frac{f_n(x+1+a)}{f_n(x+a) }.\]
Therefore, it is enough if
\[  \frac{f_n(x+1)}{f_n(x) }\leq  \frac{f_n(x+2)}{f_n(x+1) }\] for all $x,$
which holds from the assumption.
\end{proof}
\end{lemma}

\begin{lemma} \label{lem2} Under Assumption \ref{assumption}, 
$D_n(A,\{i\} )=L_n(A)-L_n(A\cup \{i\}) >0$ for any $A, \{i\} \subset \mathbb{N}, i\notin A.$

    \begin{proof} For notational convenience, we will replace $D_n$ and $ L_n$ with $D$ and  $L$ in the proof and the rest of the paper.
    
    If $A=\emptyset,$ $D(A,\{i\} )=1/(b f_n(\lfloor nT\rfloor))-1>0.$
    If $A=\{i'\} , D(A,\{i\} )=1-f_n( |i'-i| )   >0 $ since $f_n \in (0,1).$
    
    Let $|A|=2.$ We will show that
    $D(\{i_0,i_2\} ,\{i_1\} )=L(\{i_0,i_2\})-L(\{i_0,i_1, i_2\}) >0$ for any $i_0,i_1,i_2 \in \mathbb{N}.$ If $ \max\{i_0,i_2\}<i_1$ or $ \min\{i_0,i_2\}>i_1,$ the result easily follows since
    \[ L(\{i_0,i_1, i_2\})= \begin{dcases}
        L(\{i_0,i_2\})f_n(|i_1-\min\{i_0,i_2\}  |) & \text{if } i_1 <\min\{i_0,i_2\} ,\\
        L(\{i_0,i_2\}) f_n(|i_1-\max\{i_0,i_2\}  |) & \text{if } i_1 >\max\{i_0,i_2\}.
    \end{dcases} \]Now, assume $i_0<i_1<i_2.$ 
We need to show that \[ \frac{L(\{i_0,i_2 \})}{L(\{i_0,i_1,i_2\})}=\frac{f_n(x+a) }{f_n(x) f_n(a)  } >1,\] for any $x=i_1-i_0 ,a=i_2-i_1 \in \mathbb{N} .$ 
By Lemma \ref{lem1}, 
\[ \frac{f_n(x+a) }{f_n(x)f_n(a)  } \geq  \frac{f_n(1+a) }{f_n(1)f_n(a)  } \geq \frac{f_n(2) }{f_n(1)f_n(1)  }>1  .\] 

For $|A|>2, $ it is derived in the same way as $|A|=2.$ 
    If $ \max A<i$ or $ \min A>i,$ 
    \[ L(A\cup \{i\})= \begin{dcases}
        L(A)f_n(|i-\min A  |) & \text{if } i <\min A ,\\
        L(A) f_n(|i-\max A  |) & \text{if } i >\max A,
    \end{dcases} \] therefore, $D(A, \{i\}) >0.$
If $\min A <i < \max A,$ let $i^*=\min\{j: j\in A, j>i \} $, $i_*=\max\{j: j\in A, j<i \}. $
Then, \[\frac{L(A)}{L(A\cup \{i\}) }= \frac{L(\{ i_*, i^* \})}{L(\{i_*, i, i^* \}) }>1. \] 

    \end{proof}
\end{lemma}



{\bf {Proof of Theorem \ref{prop0}.}}\\
The proof is similar to the proofs of Proposition 2.2 in \cite{lee_2021_generalized} and  Theorem 3.1 in \cite{lee_2024_on}.
Note that the expression in (\ref{p2}) can be written as 
\[P ( \cap_{i'\in F }\{\xi_{i'}^{(n)}=0\}  \cap_{i\in B}   \{\xi_{i}^{(n)}=1\} ) = b f_n(\lfloor nT\rfloor) D(B,F).
\]
Therefore, if  $D(B,F)>0$ for any disjoint sets $B\sqcup F\subset [\lfloor nT\rfloor],$ the joint probability distributions are well defined, and by  Kolmogorov existence theorem, there exists a sequence of binary variables whose finite-dimensional distributions are defined as (\ref{p2}).
We will show
\begin{equation} \label{eqnD}
    D(B, F)>0
\end{equation}by mathematical induction. 
If $F=\emptyset, D(B,F)=L(B)>0.$ Assume $F\neq \emptyset.$

If $|F|=1,$ $ D(B,F)=L(B) - L(B\cup F)>0,$ by Lemma \ref{lem2}. 
Assume (\ref{eqnD}) holds with $|F|\leq m.$
We will show that (\ref{eqnD}) holds when $|F|=m+1.$ Let $F=\{i_0', i_1', i_2' , \cdots, i_m'\}$ with $i_0'<i_1'<i_2'<\cdots<i_m'.$ By definition of $D$, it follows that
\begin{align*}
D(B,F)&=  \sum_{j=0}^{|F|-1}\sum_{ \substack{F'\subset F\setminus \{i_m'\}\\ |F'|=j}}(-1)^{j} L(B\cup F')+\sum_{j=0}^{|F|-1}\sum_{ \substack{F''\subset F\setminus \{i_m'\}\\ |F''|=j}}(-1)^{j+1} L(B\cup F''\cup\{i_m'\})\\&=  D(B,F\setminus \{ i_m'\}) -
  D(B\cup \{i_m'\}, F\setminus\{ i_m'\}),\end{align*}  
therefore,
it is enough to show that \[  D(B,F\setminus \{ i_m'\}) >
  D(B\cup \{i_m'\}, F\setminus\{ i_m'\})  .\]
By
Lemma 5.1 in \cite{lee_2021_generalized}, 
\begin{equation} \label{eqn3}
    \frac{D(B,F\setminus\{ i_m'\})}{  D(B\cup\{i_m'\},F\setminus\{ i_m'\}) } =\frac{a_0-a_1-a_2-a_3\cdots-a_{m-1}}{ a_0'-a_1'-a_2'-a_3'\cdots-a_{m-1}'} , \end{equation}
where 
\begin{align*}
{a_{j}}= \begin{dcases} {  D(B\cup \{i_{j}'\},\{ i_0',i_1',\cdots,i_{j-1}'\}) } &\text{for } j=1,\cdots, m-1,\\ {  D(B,\{i_0'\}) } &\text{for } j=0, \end{dcases} \end{align*} and
\begin{align*}
{a_{j}'}= \begin{dcases} {   D(B\cup\{i_{j}',i_{m}'\},\{ i_0',i_1',\cdots,i_{j-1}'\}) } &\text{for } j=1,\cdots, m-1,\\ {   D(A_1\cup\{i_{m}'\},\{ i_0'\}) }, &\text{for } j=0. \end{dcases}
\end{align*}
We will apply Lemma 5.2 {\it i)} in \cite{lee_2021_generalized}  to show that $(6.2)>1.$ Note that $ a_0-a_1-\cdots -a_{j} >0$ 
and 
$ a_0'-a_1'-\cdots -a_{j}' >0$ for $j=0,1,\cdots,m-1,$ by the earlier assumption that (\ref{eqnD}) holds for $|F|\leq m.$

 For $j=1,\cdots, m-1,$
\[\frac{a_{j}}{a_{j}'}=\frac{  D(B\cup \{i_{j}'\},\{ i_0',i_1',\cdots,i_{j-1}'\}) }{  D(B\cup\{i_{j}',i_{m}'\},\{ i_0',i_1',\cdots,i_{j-1}'\}) }= \frac{\sum_{\ell=0}^{j}\sum_{ \substack{F'\subset  \{i_0',i_1',\cdots,i_{j-1}'\}\\ |F'|=\ell}}(-1)^{\ell} L(B\cup\{i_j'\}\cup F')}{\sum_{\ell'=0}^{j}\sum_{ \substack{F''\subset \{ i_0',i_1',\cdots,i_{j-1}' \}\\ |F''|=\ell'}}(-1)^{\ell'} L(B\cup\{i_j', i_m'\} \cup F'')} .\]
Note that
if $B=\emptyset,$ then, 
for any $F'=F''\subset  \{i_0',i_1',\cdots,i_{j-1}'\},$
\[ \frac{L(B\cup\{i_j'\} \cup F')}{  L(B\cup\{i_j', i_m'\} \cup F'')}  = \frac{1}{f_n(i_{m}'-i_j') }.\]
Similarly, for any $F'=F''\subset  \{i_0',i_1',\cdots,i_{j-1}'\},$ if  $\max B <i_m'$, define  $i_j^*= \max \{ i_j', i: i\in B, i<i_m' \},$ then
\[ \frac{L(B\cup\{i_j'\} \cup F')}{  L(B\cup\{i_j', i_m'\} \cup F'')}  = \frac{1}{f_n(i_{m}'-i_j^*) }.\] 
If $\max  B >i_m',$ define $i^*= \min \{i: i\in B, i>i_m' \},$ then
\[ \frac{L(B\cup\{i_j'\} \cup F')}{  L(B\cup\{i_j', i_m'\} \cup F'')}  = \frac{f_n(i^*-i_j^*)}{f_n(i_{m}'-i_j^*) f_n(i^*-i_m' ) } .\]

Therefore,
\begin{align} \label{eqnD2}
\frac{a_{j}}{a_{j}'} =\begin{dcases}\frac{f_n(i^*-i_j^*)}{f_n(i_{m}'-i_j^*)f_n(i^*-i_m' ) } &\text{ if }  \max  B >i_m',  \\ \frac{1}{f_n(i_{m}'-i_j^*) } & \text{ if }   \max B <i_m',\\
\frac{1}{f_n(i_{m}'-i_j') } & \text{ if }   B=\emptyset.
\end{dcases}  \end{align}  
Also,
\begin{equation*}
  \frac{a_{0}}{a_{0}'}=\frac{  D(B,\{i_0'\}) }{   D(B\cup\{i_{m}'\},\{ i_0'\}) }=\frac{L(B)-L(B\cup \{i_0'\})}{ L(B\cup \{i_m'\})- L(B\cup \{i_0', i_m'\})  }
\end{equation*}
with
\begin{equation*}
    \frac{L(B)}{L(B\cup \{i_m'\})}=\begin{dcases}
    \frac{1}{f_n(i^*-i_{m}') } &\text{ if } \min B>i_m',
    \\ \frac{1}{f_n(i_{m}'-i_*) } &\text{ if } \max B<i_m' ,\\
  \frac{f_n(i^*-i_*)}{f_n(i^*-i_{m}')f_n(i_m'-i_*) } &\text{ if } \min B< i_m'<\max B, \\
   \frac{1}{b f_n(\lfloor nT\rfloor)} &\text{ if } B=\emptyset,
    \end{dcases}
\end{equation*} where $i_*=\max\{i: i\in B, i<i_m'\}$ if $\min B<i_m',$ and 

\begin{equation*}
 \frac{L(B\cup \{i_0'\})}{ L(B\cup \{i_0', i_m'\})  }=(\ref{eqnD2})
\end{equation*} with $j=0.$
Since $i_*\leq i_0^*\leq i_1^* \leq i_2 ^* \leq \cdots \leq i_{m-1}^*< i_m'<i^*,$ when they are defined, and $f_n(x)$ is a decreasing function such that  $f_n(x+a)/f_n(x) $ is  a non-decreasing function of $x$ for any $a\in \mathbb{N}$ by Lemma \ref{lem1}, and  $f_n(x+a)/f_n(x)<1 $ for any $x,a \in \mathbb{N}$,  it follows that

\begin{align*}
  &\frac{L(B)}{ L(B\cup \{i_m'\})}  \geq \frac{L(B\cup \{i_0'\})}{ L(B\cup \{i_0', i_m'\})  } \geq  \frac{a_1}{a_1'}  \geq  \frac{a_2}{a_2'} \geq \cdots   \geq  \frac{a_{m-1}}{a_{m-1}'}
\end{align*}
therefore, by Lemma 5.2  {\it i)} in \cite{lee_2021_generalized} and Lemma \ref{lem2}, 
\[(\ref{eqn3})\geq  \frac{L(B)}{ L(B\cup \{i_m'\})}>1. \]
\qed \\ \\

We will denote the discrete convolution as follows. \[ (f_n*f_n)[j]= \sum_{\ell=1}^{j-1} f_n(j-\ell)f_n(\ell),\] and for $k\geq 2,$
$ f_n^{*k}[j]=  f_n*(f_n^{*(k-1)}) [j] .$
 Under Assumption \ref{real_assum}, 
$f_n(\lfloor nT \rfloor)\sim ce^{-\lambda T}\lfloor nT \rfloor^{\alpha-1},$  for any $T>0,$ as $n\to \infty$.   Also, 
for any $k\in \mathbb{N}, t>0,$ as $ n \to \infty$,
    \begin{equation}\label{approx0}
    f_n^{*k}[tn]\sim \frac{ e^{-\lambda t}(c\Gamma(\alpha)) ^k}{\Gamma(k\alpha )} \lfloor tn\rfloor ^{k\alpha -1}.\end{equation}
    
{\bf {Proof of Theorem \ref{theo1}.}}\\
First, we will prove that $X_t^{(n)}$ converges in distribution and its limiting distribution has mgf (\ref{mgf1}). Similarly, we will show that all finite-dimensional distributions of  $\{ X_t^{(n)}, 0<t<T\}$ converge weakly and its limiting distribution has mgf (\ref{mgf2}).
Then, by Kolmogorov's existence theorem, there exists a stochastic process $X_t$ whose finite-dimensional distributions are defined by (\ref{mgf2}). We will show that $X_t$ has stationary increments and is continuous in probability.
After showing that \\
{\it i)} all finite-dimensional distributions of $X_t^{(n)}$ converge to the corresponding finite-dimensional distributions of $X_t$,\\
{ \it ii)} $X_t$ is continuous in probability, and\\
{\it iii)} $X_t^{(n)}$ has non-decreasing sample path,\\ 
it follows that $\{ X_t^{(n)}, 0<t<T\}$ converges weakly to $\{X_t, 0<t<T\}$  in the space of $D[0,T]$ with Skorokhod ($J_1$) topology by Theorem 3 in \cite{bingham_1971_limit}.

Step 1) To show the convergence of the marginal distribution.\\
Note that \begin{align*}
& n^{\alpha k}\mathrm{E}[ (X_t^{(n)})^k  ] = \mathrm{E}[(  \sum_{i\in \lfloor nt \rfloor  } \xi_i^{(n)}  )^k]=    k! \sum_{1\leq i_1<i_2<\cdots<i_k\leq  \lfloor nt \rfloor  } \mathrm{E}\Big[\prod_{j=1}^{k}\xi_{i_j}^{(n)}  \Big]\\&+  \sum_{\substack{\ell_i\geq 1\\ \ell_1+\ell_2+\cdots+\ell_{k-1}=k}}\frac{k!}{\ell_1!\ell_2!\cdots \ell_{k-1}!} \sum_{1\leq i_1<i_2<\cdots<i_{k-1}\leq  \lfloor nt \rfloor  } \mathrm{E}\Big[\prod_{j=1}^{k-1}\xi_{i_j}^{(n)}  \Big]\\&+  \sum_{\substack{\ell_i\geq 1\\ \ell_1+\ell_2+\cdots+\ell_{k-2}=k}}\frac{k!}{\ell_1!\ell_2!\cdots \ell_{k-2}!} \sum_{1\leq i_1<i_2<\cdots<i_{k-2}\leq  \lfloor nt \rfloor  } \mathrm{E}\Big[\prod_{j=1}^{k-2}\xi_{i_j}^{(n)}  \Big] \\& +\cdots +
  \sum_{1\leq i\leq  \lfloor nt \rfloor  }    \mathrm{E} [ \xi_{i}^{(n)}  ] \end{align*}  
  \begin{align*}
    & =k!\sum_{\substack{1\leq  i_1<i_2<\cdots<i_k \leq \lfloor nt \rfloor } } P(\cap_{j=1,\cdots, k} \{X_{i_j}=1\})\\&+\sum_{\substack{\ell_i\geq 1\\ \ell_1+\ell_2+\cdots+\ell_{k-1}=k}}\frac{k!}{\ell_1!\ell_2!\cdots \ell_{k-1}!}\sum_{\substack{i_1,\cdots,i_{k-1}\in \lfloor nt \rfloor\\ i_1<i_2<\cdots<i_{k-1} } } P(\cap_{j=1,\cdots, k-1} \{X_{i_j}=1\})\\&+\sum_{\substack{\ell_i\geq 1\\ \ell_1+\ell_2+\cdots+\ell_{k-2}=k}}\frac{k!}{\ell_1!\ell_2!\cdots \ell_{k-2}!}\sum_{\substack{i_1,\cdots,i_{k-2}\in \lfloor nt \rfloor \\ i_1<i_2<\cdots<i_{k-2} } } P(\cap_{j=1,\cdots, k-2} \{X_{i_j}=1\})\\&+\cdots +\sum_{i=1}^{\lfloor nt \rfloor} P( \{X_{i}=1\}) 
   \\&
     = k!\sum_{\ell=k-1}^{\lfloor nt\rfloor -1}  \sum_{j=k-1}^{\ell} b f_n(\lfloor nT \rfloor) f_n^{*(k-1)}[j] \\&
+\sum_{\substack{\ell_i\geq 1\\ \ell_1+\ell_2+\cdots+\ell_{k-1}=k}}\frac{k!}{\ell_1!\ell_2!\cdots \ell_{k-1}!} \sum_{\ell=k-2}^{\lfloor nt\rfloor -1}  \sum_{j=k-2}^{\ell}  b f_n(\lfloor nT \rfloor) f_n^{*(k-2)}[j]\\&+\cdots + \lfloor nt\rfloor b f_n(\lfloor nT \rfloor)\\&= \lfloor nt\rfloor b f_n(\lfloor nT \rfloor)  +\sum_{m=1}^{k-1} (m+1)! S(k,m+1) \sum_{\ell=m}^{\lfloor nt\rfloor -1} \sum_{j=m}^{\ell } b f_n(nT) f_n^{*m}[j], \end{align*} 
where $S(\cdot, \cdot)$ denotes the Stirling number of the second kind.

Under Assumption \ref{real_assum}, and by (\ref{approx0}), we have
\[a_1= \lim_{n\to \infty}\lfloor nt\rfloor \frac{b f_n(\lfloor nT \rfloor)}{n^{\alpha}}= cb t e^{-\lambda T}T^{\alpha-1}, \]

and  \[a_k=\lim_{n\to \infty}\sum_{\ell=k-1}^{\lfloor nt\rfloor -1} \sum_{j=k-1}^{\ell } b f_n(nT) \frac{f_n^{*(k-1)}[j]}{n^{\alpha k}}= b c   e^{-\lambda T}T^{\alpha-1}
 \int_{0}^t \int_{0}^x  \nu^{k-1}  \frac{e^{-\lambda y}y^{(k-1)\alpha-1}}{\Gamma((k-1)\alpha)} dy    dx. \]
Therefore, 
  $ \lim_{n\to \infty} \mathrm{E}[ (X_t^{(n)})^k  ]=k!a_{k}$ for all $k\in \mathbb{N},$
    $\{ X_t^{(n)}, n\in \mathbb{N} \} $ is tight, and every subsequence of $ X_t^{(n)}$ converges to a random variable  whose k-th moment is $k! a_k$    \citep{durrett_2010_probability}. The convergence of $ \sum_{k\in \mathbb{N}} s^k a_k$ for any $s\in \mathbb{R}$ can be shown using the Mittag-Leffler function, which leads to a unique limiting distribution. 
    Therefore, $X_t^{(n)}$ converges to the random variable $X_t$ 
    whose mgf is
\[  \mathrm{E}( e^{s X_t})= \sum_{k=0}^{\infty} a_ks^k, \] where $a_0=1.$

Step 2) To show the convergence of all finite-dimensional distributions.\\
Assume that
for any $m\in \mathbb{N}, C\subset [m], C\neq \emptyset,$ and $ k_i\in \mathbb{N} , i\in [|C|],$
\begin{align} \label{equ4}
  &\sum_{\substack{ \lfloor nt_{C(1)-1} \rfloor <i_1^1<i_2^1<\cdots <i_{k_1}^1\leq \lfloor nt_{C(1)} \rfloor \\
     \lfloor nt_{C(2)-1}\rfloor<i_{1}^2<i_{2}^2<\cdots <i_{k_2}^2\leq  \lfloor nt_{C(2)} \rfloor \\ \vdots\\
  \lfloor nt_{C(|C|)-1} \rfloor <i_{1}^{|C|}<\cdots <i_{k_{|C|} }^{|C|}\leq \lfloor nt_{C(|C|)} \rfloor  } }  \mathrm{E}\Big(  \prod_{j\in [|C|]}\prod_{\ell \in [ k_j] } 
 \xi_{i_{\ell}^{j}  }^{(n)} 
\Big) \nonumber \\& \sim  n^{\alpha (k_1+k_2+\cdots +k_{|C|})} a_T   t_{C(|C|)}^{\alpha(\sum_{i\in  [|C|] }k_i-1)+1 }  A(C,\{k_i\})  ,\end{align}
then we have
\begin{align*}
  & \lim_{n\to \infty}\mathrm{E}\Bigg( \prod_{j\in [|C|]} \frac{ (X_{t_{C(j)}}^{(n)}-X_{t_{C(j)-1} }^{(n)})^{k_j}}{k_j!}   \Bigg)\\&=
\lim_{n\to \infty}  \sum_{\substack{k_j'\in [k_j]\\ j\in [|C|]}} \Big(\prod_{j\in [|C|]}\frac{k_j'!S(k_j,k_j')}{k_j!} \Big)
    \sum_{ \substack{\lfloor nt_{C(j)-1} \rfloor < i_1^j<\cdots<i_{k_j'}^j \leq \lfloor nt_{C(j)} \rfloor\\ j\in [ |C| ]}   }  \mathrm{E}\Big(  \prod_{j\in [|C|]}\prod_{\ell \in [ k_j'] } 
 \xi_{i_{\ell}^{j}  } ^{(n)}
\Big) \frac{1}{ n^{\alpha (k_1+k_2+\cdots +k_{|C|})}}
  \\&=\lim_{n\to \infty} 
    \sum_{ \substack{\lfloor nt_{C(j)-1} \rfloor < i_1^j<\cdots<i_{k_j}^j \leq \lfloor nt_{C(j)} \rfloor\\ j\in [ |C| ]}   }  \mathrm{E}\Big(  \prod_{j\in [|C|]}\prod_{\ell \in [ k_j] } 
 \xi_{i_{\ell}^{j}  } ^{(n)}
\Big) \frac{1}{ n^{\alpha (k_1+k_2+\cdots +k_{|C|})}}\\ &= a_T   t_{C(|C|)}^{\alpha(\sum_{i\in  [|C|] }k_i-1)+1 }  A(C,\{k_i\}), \end{align*} where the last two equalities are derived by  (\ref{equ4}), and the tightness of all finite-dimensional distributions of $X_t^{(n)}$ follows.  If \[ \sum_{\substack{k_i\in\mathbb{N}\\ i\in [ |C| ]}} t_{C(|C|)}^{\alpha(\sum_{i\in [|C|]}k_i-1)+1 } \Big(\prod_{i\in [|C| ] }S^{k_i}\Big) A( C, \{k_i\} )<\infty\] for any $C\subset \mathbb{N}, S\in \mathbb{R}_+,$
then the limiting distribution of all finite-dimensional distributions of $X_t^{(n)}$  is uniquely defined with mgf $(\ref{mgf2}).$
By the change of variables, $A(C,\{k_i\})$ can be expressed as convolutions of beta distributions as follows.
For simplicity, letting $ k_i>1, i\in [|C|],$ we have
\begin{align*}
    A(C,\{k_i\})&=\int^1_{t_{C(|C|)-1}/t_{C(|C|)}} \int_{t_{C(|C|)-1}/t_{C(|C|)} }^{x_{|C|}} \cdots \cdots 
\int_{t_{C(2)-1}/t_{C(|C|)}}^{t_{C(2)}/t_{C(|C|)}}
    \int_{t_{C(2)-1}/t_{C(|C|)}}^{x_2}
\int_{t_{C(1)-1}/t_{C(|C|)}}^{t_{C(1)}/t_{C(|C|)}} \int_{t_{C(1)-1}/t_{C(|C|)}}^{x_1} \\&
\frac{\nu^{k_1-1}(x_1-x_1')^{\alpha(k_1-1)-1}}{\Gamma((k_1-1)\alpha ) } \prod_{i=2}^{|C|} 
 \Big\{ c ( x_{i}'-x_{i-1})^{\alpha-1} \frac{ \nu^{k_i-1} (x_{i}-x_{i}')^{\al(k_i-1)-1}}{\Gamma((k_i-1)\al)}\Big\}\\& \exp(-\lambda t_{C(|C|)} ( x_{|C|}-x_1'))
    dx_1'dx_1dx_2'dx_2\cdots dx_{|C|}' dx_{|C|}.\end{align*}
    Using the change of variables with 
    $y_i'=x_i'/x_i, y_i=x_i/x_{i+1}'$ for $i<|C| ,$ and $y_{|C|}=x_{|C|},$ 
    \begin{align*}
    A(C,\{k_i\})&=
   \frac{\nu^{   \sum_{i\in [|C|]}k_i-1}}{\Gamma\big(\alpha (\sum_{i\in [|C|] }k_i-1)+1\big)}\iiint\cdots \int_{Q_C^y}  
    f(y_1',1,\al(k_1-1))f(y_1, \al(k_1-1)+1,\al) 
    \\&\times \prod_{2\leq i<|C|}\Big\{f\Big(y_i', \al(\sum_{j=1}^{i-1} k_j)+1, \al(k_i-1)\Big) f\Big(y_i, \al(\sum_{j=1}^i k_j-1)+1, \al\Big)\Big\} \\& \times  f\Big(y_{|C|}', \al(\sum_{j=1}^{|C|-1} k_j)+1, \al(k_{|C|}-1)\Big) \\&\times \exp\Big\{-\lambda t_{C(|C|)} \Big(y_{|C|}'-  \prod_{i=1}^{|C|}y_iy_i'\Big)\Big\} 
    dy_1'dy_1dy_2'dy_2\cdots dy_{|C|}' dy_{|C|},
\end{align*}
where  $f(x, a,b )$ is the pdf of the beta distribution,
\[
f(x, a,b )=\frac{\Gamma(a+b)x^{a-1}(1-x)^{b-1}}{ \Gamma(a)\Gamma(b)},\]
and 
\[Q_C^y=\Big\{ y_i,y_i': y_i'=\frac{x_i'}{x_i}, y_i=\frac{x_i}{x_{i+1}'},   x_{|C|+1}'=x_{|C|},\frac{t_{C(i)-1}}{t_{C(|C|)}}<x_i'<x_i< \frac{t_{C(i)}}{t_{C(|C|)}}  \text{ for } i\in[|C|]  \big\}.  \]
Therefore, \begin{align} \label{A}
    A(C,\{k_i\})\leq \frac{\nu^{   \sum_{i\in [|C|]}k_i-1}}{\Gamma\big(\alpha (\sum_{i\in [|C|] }k_i-1)+1\big)} (e^{-\lambda t_{C(|C|)}} \vee 1),
\end{align}   and one can show in a similar way that (\ref{A}) holds for any $\{k_i \in \mathbb{N}\}, C\subset\mathbb{N}.$ Using the three-parameter generalization of Mittag-Leffler function, if follows that for any $C\subset [m], S\in \mathbb{R}_+,$ 
\[ \sum_{\substack{k_i\in \mathbb{N}\\ i\in [|C|]} }    (St_{C(|C|)}^{\alpha})^{(\sum_{i\in  [|C|] }k_i-1) } \frac{\nu^{   \sum_{i\in [|C|]}k_i-1}}{\Gamma\big(\alpha (\sum_{i\in [|C|] }k_i-1)+1\big)}  <(St^{\al}_{C(|C|)}\nu)^{-1} E_{\al,1-\al}^{|C|}( St_{C(|C|)}^{\al}\nu)<\infty, \] where $ E_{\al,\be}^{\gamma}(z)$ is the three-parameter Mittag-Leffler function,
\[ E_{\al,\be}^{\gamma}(z)=  \frac{1}{\Gamma(\gamma) } \sum_{k=0}^{\infty} \frac{\Gamma(\gamma+k)z^k}{k! \Gamma(\al k+\be)}.\]

 Therefore, $(\ref{mgf2})$ converges,
and the weak convergence of all finite-dimensional distributions of $X_t^{(n)}$ follows with the mgf of a finite-dimensional limiting distribution as (\ref{mgf2}).

Now we will show that (\ref{equ4}) holds. Assume $k_1, k_2,\cdots,k_{|C|}>1.$ Other cases are proved similarly. Note that
\begin{align} \label{equ5}
    & \sum_{\substack{\lfloor nt_{C(j)-1} \rfloor < i_1^j<\cdots<i_{k_j}^j \leq \lfloor nt_{C(j)} \rfloor\\ j\in [ |C| ]}} \mathrm{E}\Big(  \prod_{j\in [|C|]}\prod_{\ell \in [ k_j] } 
 \xi_{i_{\ell}^{j}  } ^{(n)}
\Big)  \nonumber
\\ &=\sum_{\substack{\lfloor nt_{C(j)-1} \rfloor < i_1^j<i_{k_j}^j \leq \lfloor nt_{C(j)} \rfloor\\ j\in [ |C| ]}}  b f_n( \lfloor nT \rfloor)f_n^{*k_1-1}[i_{k_1}^1-i_1^1]f_n[i_1^2-i_{k_1}^1] f_n^{*k_2-1}[i_{k_2}^2-i_1^2]f_n[i_1^3-i^2_{k_2}]\cdots \nonumber \\& \hspace{95pt}\times f_n[i_{{1}}^{|C|}-i_{{k_{|C|-1}}}^{|C|-1}]f_n^{* k_{|C|} -1}[i_{k_{|C|}}^{|C|} -i_{{1}}^{|C|}]  ,
\end{align}
where $t_0=0.$ Using (\ref{approx0}) and the change of variable,
\begin{align*}
    (\ref{equ5})\sim &\sum_{\substack{\lfloor nt_{C(j)-1} \rfloor < i_1^j<i_{k_j}^j \leq \lfloor nt_{C(j)} \rfloor\\ j\in [ |C| ]}} \frac{ \nu^{k_1-1}}{\Gamma(\alpha(k_1-1 ))} \lfloor i_{k_1}^1-i_1^1\rfloor ^{\alpha(k_1-1) -1} \\&\times \prod_{j=2}^{|C|} \Bigg\{ c  \lfloor i_{1}^j-i_{k_{j-1}}^{j-1}\rfloor ^{\alpha -1}  \frac{ \nu^{k_j-1}}{\Gamma(\alpha(k_j-1 ))} \lfloor i_{k_j}^j-i_1^j\rfloor ^{\alpha(k_j-1) -1}\Bigg\}\exp\big(-\lambda (i_{k_{|C|}}^{|C|} -i_1^1)/n\big)
    \\ &\sim n^{\alpha (k_1+k_2+\cdots +k_{|C|})} a_T   t_{C(|C|)}^{\alpha(\sum_{i\in  [|C|] }k_i-1)+1 }   \iint \cdots \int_{Q_C}
 \frac{\nu^{k_1-1}(x_1-x_1')^{\alpha(k_1-1)-1}}{\Gamma((k_1-1)\alpha ) }   \\ & \times \prod_{i=2}^{|C|} 
 \Big\{ c ( x_{i}'-x_{i-1})^{\alpha-1} \frac{ \nu^{k_i-1} (x_{i}-x_{i}')^{\al(k_i-1)-1}}{\Gamma((k_i-1)\al)}\Big\}
  \exp\big({-\lambda t_{C(|C|)} (x_{|C|}-x_1')}\big)
   \prod_{\ell\in [|C|]}dx_{\ell}'dx_{\ell}  \\ &=n^{\alpha (k_1+k_2+\cdots +k_{|C|})} a_T   t_{C(|C|)}^{\alpha(\sum_{i\in  [|C|] }k_i-1)+1 }  A(C,\{k_i\}),
\end{align*}
where
\begin{align*}
 Q_C=\Big\{ x_i',x_i, x_j:  
\frac{ t_{C(i)-1}}{t_{C(|C|)}}
 <x_i'<x_i<\frac{t_{C(i)}}{t_{C(|C|)}}, i\in  [|C|] 
 \Big\}. \end{align*}

Step 3) By Kolomogorov existence theorem, there exists a stochastic process $X_t$ such that all finite-dimensional distributions are specified as (\ref{mgf2}). 
Since $X_t^{(n)}$ has stationary increments and its finite-dimensional distributions converge to those of $X_t,$ it follows that $X_t$ has stationary increments. Also, from (\ref{mgf1}), $\mathrm{E}(X_t^2)\to0$ as $t\to 0$, from which it follows that for any $\epsilon>0,$
\[   P(  |  X_t-X_{t_n}  |>\epsilon) \to 0 \] as $t_n\to t$, therefore, $X_t$ is continuous in probability. It is trivial that $X_t^{(n)}$ has a non-decreasing sample path. By Theorem 3 in \cite{bingham_1971_limit}, 
\[   X_t^{(n)}=\frac{\sum_{i\in \lfloor nt \rfloor  } \xi_i^{(n)}  }{ n^{\alpha} } \Rightarrow X_t\]  in the space $D[0,T]$ with the Skorokhod ($J_1$)  topology.
\qed \\

{\bf {Proof of Corollary \ref{cor1}.}}\\
  {\it i,ii)} The results can easily be derived by (\ref{mgf1}) and (\ref{mgf2}).
 \\ {\it iii)} Using (\ref{mgf2})
    with $ C=\{ 2,4\}, k_1=k_2=1,$ one can show that
\begin{align*}
 \mathrm{E}((X_{t_2}-X_{t_1})(X_{t_4}-X_{t_3}))=a_T t_{4}^{\alpha+1} \int_{t_3/t_4}^{1}\int_{t_1/t_4}^{t_2/t_4 } c(x_2-x_1)^{\al-1}e^{-\lambda t_4(x_2-x_1) } dx_1 dx_{2},\end{align*}
from which we obtain, for $d\approx 0,$ 
\[ \mathrm{E}(\nabla X_{s}  \nabla X_{t+s} )\approx (a_Tc) d^2  e^{-\lambda t} t^{\alpha-1} . \] 
Also, using (\ref{mgf1}), it follows that
\[\mathrm{E}(\nabla X_{s}^k)\approx d^2 \frac{k!a_T\nu^{k-1}}{\Gamma((k-1)\alpha+2)}  d^{(k-1)\alpha-1}, \] and the results follow.
\\{\it iv)} For simplicity, let $k_i>1, i=1,2,\cdots,m.$ By (\ref{mgf2}) and using the change of variable, 

\begin{align*}
  \mathrm{E}\Big(\prod_{i=1}^m (\nabla X_{t_i})^{k_i}\Big) =a_T\Big(\prod_{i=1}^mk_i !\Big)&\int_{t_m-d}^{t_m} \int_{t_m-d}^{x_m} \cdots   \int_{t_2-d}^{t_2} \int_{t_2-d}^{x_2} \int_{t_1-d}^{t_1} \int_{t_1-d}^{t_1}
  \frac{\nu^{k_1-1}(x_1-x_1')^{\alpha(k_1-1)-1}}{\Gamma((k_1-1)\alpha ) } \\& \times \prod_{i=2}^m 
 \Big\{ c ( x_{i}'-x_{i-1})^{\alpha-1} \frac{ \nu^{k_i-1} (x_{i}-x_{i}')^{\al(k_i-1)-1}}{\Gamma((k_i-1)\al)}\Big\}\\ &\times
\exp\big({-\lambda  (x_{m}-x_1')}\big) dx_1'dx_1 dx_2'dx_2  \cdots dx_m' dx_m.
\end{align*}
As $d\to 0,$ it is approximated as 
\begin{align}\label{approx}
  \mathrm{E}\Big(\prod_{i=1}^m (\nabla X_{t_i})^{k_i}\Big) \sim &a_T\Big(\prod_{i=1}^mk_i !\Big) 
  \frac{\nu^{k_1-1} d^{\alpha (k_1-1)+1}}{\Gamma((k_1-1)\alpha+2 ) }  \prod_{j=2}^m \Big\{c(t_j-t_{j-1} )^{\al-1} \frac{\nu^{k_j-1} d^{\alpha (k_j-1)+1 }}{\Gamma((k_j-1)\alpha+2 ) } \Big\}\\ &\times \exp(-\lambda (t_m-t_1)) . \nonumber
\end{align}
One can show in a similar way that (\ref{approx}) holds for any $k_i\in \mathbb{N}, i=1,2,\cdots,m.$
Therefore, the results follows.
    \\ \qed \\

{\bf {Proof of Proposition \ref{prop_01}.}}\\
Note that the expression in (\ref{p3}, \ref{p4}) can be written as 
\begin{align*} P ( \cap_{i'\in F }\{\xi_{i'}^{(*,n)}=0\}  \cap_{i\in B}   \{\xi_{i}^{(*,n)}=1\} ) =  D(B_{+1},F_{+1}'),
\end{align*}  where $B_{+1}:=\{1, i+1 : i\in B \}$ and $F_{+1}':=\{ i+1: i\in F\}.$  If $B=\emptyset, B_{+1}:=\{1\}, $ and if $F=\emptyset, F_{+1}':= \emptyset.$ Note that \begin{align*}
     \sum_{ B\sqcup F\subset [n]  } D(B_{+1},F_{+1}') =1\end{align*}  for any $n\in \mathbb{N}$. If  $ D(B_{+1},F_{+1}')  >0$ for any disjoint sets $B\sqcup F\subset \mathbb{N},$ the  probability distributions  are well defined, and by  Kolmogorov existence theorem, there exists a sequence of binary variables whose finite-dimensional distributions are defined as (\ref{p3}, \ref{p4}). 
Following the proof of Theorem \ref{prop0}, one can easily show that $D(B_{+1},F_{+1}')>0$ for any disjoint sets $B\sqcup F\subset \mathbb{N}.$  
\qed \\

{\bf {Proof of Theorem \ref{theo2}.}}\\
It will be proved in the same way as Theorem \ref{theo1}.\\
Step 1) To show the weak convergence of marginal distribution.

Note that \begin{align*}
     n^{\alpha k}\mathrm{E}[ (X_t^{(*,n)})^k  ] =& \mathrm{E}[(  \sum_{i\in \lfloor nt \rfloor  } \xi_i^{(*,n)}  )^k]
     = k! \sum_{1\leq i_1<i_2<\cdots<i_k\leq  \lfloor nt \rfloor  } \mathrm{E}\Big[\prod_{j=1}^{k}\xi_{i_j}^{(*,n)}  \Big]\\&+  \sum_{\substack{\ell_i\geq 1\\ \ell_1+\ell_2+\cdots+\ell_{k-1}=k}}\frac{k!}{\ell_1!\ell_2!\cdots \ell_{k-1}!} \sum_{1\leq i_1<i_2<\cdots<i_{k-1}\leq  \lfloor nt \rfloor  } \mathrm{E}\Big[\prod_{j=1}^{k-1}\xi_{i_j}^{(*,n)}  \Big]\\&+  \sum_{\substack{\ell_i\geq 1\\ \ell_1+\ell_2+\cdots+\ell_{k-2}=k}}\frac{k!}{\ell_1!\ell_2!\cdots \ell_{k-2}!} \sum_{1\leq i_1<i_2<\cdots<i_{k-2}\leq  \lfloor nt \rfloor  } \mathrm{E}\Big[\prod_{j=1}^{k-2}\xi_{i_j}^{(*,n)}  \Big] \\& +\cdots +
  \sum_{1\leq i\leq  \lfloor nt \rfloor  }    \mathrm{E} [ \xi_{i}^{(*,n)}  ] 
     \\&= \sum_{m=1}^k  m!S(k,m) \sum_{j=m}^{\lfloor nt\rfloor } f_n^{*m}[j]. \end{align*} 
With Assumption \ref{real_assum} and (\ref{approx0}),  \[a_k^*=\lim_{n\to \infty}\sum_{j=k}^{\lfloor nt \rfloor } \frac{f_n^{*k}[j]}{n^{\alpha k}}= \sum_{k=1}^{\infty}(cs\Gamma(\alpha))^{k} \int_{0}^t \frac{e^{-\lambda y}y^{k\alpha-1}}{\Gamma(k\alpha)} dy. \]
Therefore,  $ \lim_{n\to \infty} \mathrm{E}[ (X_t^{(*,n)})^k  ]=k!a_k^*$ for all $k\in \mathbb{N},$
    and $\{ X_t^{(*,n)}, n\in \mathbb{N} \} $ is tight. Every subsequence of $ X_t^{(*,n)}$ converges to a random variable  whose $k$-th moment is $k! a_k^*,$ and   it can be shown that $ \sum_{k\in \mathbb{N}} s^k a_k^*<\infty$ using Mittag-Leffler function. Therefore, the limiting distribution is uniquely defined, and $ X_t^{(*,n)}$ converges weakly to $X_t^*$ whose  mgf is
\[  \mathrm{E}( e^{s X_t^*})= \sum_{k=0}^{\infty} a_k^*s^k , \] where $a_0^*=1.$

    Step 2) To show the weak convergence
of all finite-dimensional distributions of $X_t^{(*,n)}.$ 

Assume that
for any $m\in \mathbb{N} ,C\subset [m], C\neq \emptyset,$ and $ \{k_i\in \mathbb{N} ; i\in [|C|]\},$
\begin{align} \label{equ7} &
  \sum_{  \substack{\lfloor nt_{C(j)-1} \rfloor < i_1^j<\cdots<i_{k_j}^j \leq \lfloor nt_{C(j)} \rfloor\\ j\in [ |C| ]}  }  \mathrm{E}\Big(  \prod_{j\in [|C|]}\prod_{\ell \in [ k_j] } 
 \xi_{i_{\ell}^{j}}^{(*,n)}   
\Big) \nonumber \\& \sim  n^{\alpha (k_1+k_2+\cdots +k_{|C|})}    t_{C(|C|)}^{\alpha\sum_{i\in  [|C|] }k_i } A^*( C, \{k_i\} ) ,\end{align} then we have
\begin{align*}
  & \mathrm{E}\Big( \prod_{j\in [|C|]} \frac{ (X_{t_{C(j)}}^{(*,n)}-X_{t_{C(j)-1}}^{(*,n)})^{k_j}}{k_j!}   \Big)\\&=\lim_{n\to \infty} 
    \sum_{\substack{\lfloor nt_{C(j)-1} \rfloor < i_1^j<\cdots<i_{k_j}^j \leq \lfloor nt_{C(j)} \rfloor\\ j\in [ |C| ]} }  \mathrm{E}\Big(  \prod_{j\in [|C|]}\prod_{\ell \in [ k_j] } 
 \xi_{i_{\ell}^{j}  }^{(*,n)} 
\Big) \frac{1}{ n^{\alpha (k_1+k_2+\cdots +k_{|C|})}}\\ &=  t_{C(|C|)}^{\alpha\sum_{i\in  [|C|] }k_i } A^*( C, \{k_i\} ),\end{align*} 
and the tightness of all finite-dimensional distributions of $X_t^{(*,n)}$ follows. Using a similar method used to derive (\ref{A}),  one can show that $A^*( C, \{k_i\} )$ is expressed as convolutions of beta distributions and \[
A^*( C, \{k_i\} ) \leq  \frac{\nu^{   \sum_{i\in [|C|]}k_i}}{\Gamma\big(\alpha \sum_{i\in [|C|] }k_i\big)}  (e^{-\lambda t_{C(|C|)}} \vee 1).
\]
Using the three-parameter generalization of Mittag-Leffler function, one can show that for any $C\subset [m], S\in \mathbb{R}_+,$ 
\[ \sum_{\substack{k_i\in \mathbb{N}\\ i\in [|C|]} } (St_{C(|C|)}^{\alpha})^{\sum_{i\in  [|C|] }k_i } \frac{\nu^{   \sum_{i\in [|C|]}k_i}}{\Gamma\big(\alpha \sum_{i\in [|C|] }k_i\big)}<   E_{\al,0}^{|C|}( St_{C(|C|)}^{\al}\nu)<\infty. \] 
 Therefore, (\ref{mgf4}) converges,
from which the weak convergence of all finite-dimensional distributions of $X_t^{(n)}$ follows with the mgf of all finite-dimensional limiting distributions as (\ref{mgf4}).

Next, we show that (\ref{equ7}) holds. Let $k_1, k_2,\cdots,k_{|C|}>1. $ Other cases are proved similarly.
Note that
\begin{align} \label{eqn4}
    & \sum_{\substack{\lfloor nt_{C(j)-1} \rfloor < i_1^j<\cdots<i_{k_j}^j \leq \lfloor nt_{C(j)} \rfloor\\ j\in [ |C| ]}} \mathrm{E}\Big(  \prod_{j\in [|C|]}\prod_{\ell \in [ k_j] } 
 \xi_{i_{\ell}^{j}  }^{(*,n)} 
\Big)  \nonumber
\\ &=\sum_{\substack{\lfloor nt_{C(j)-1} \rfloor < i_1^j<i_{k_j}^j \leq \lfloor nt_{C(j)} \rfloor\\ j\in [ |C| ]}}  f_n^{*k_1}[i_{k_1}^1]f_n[i_1^2-i_{k_1}^1] f_n^{*k_2-1}[i_{k_2}^2-i_1^2]f_n[i_1^3-i^2_{k_2}]\cdots \nonumber \\& \hspace{95pt} \times f_n[i_{{1}}^{|C|}-i_{{k_{|C|-1}}}^{|C|-1}]f_n^{* k_{|C|} -1}[i_{k_{|C|}}^{|C|} -i_{{1}}^{|C|}],
\end{align}
where $t_0=0.$ Using the change of variable and (\ref{approx0}), 
\begin{align*}
    (\ref{eqn4})\sim &\sum_{\substack{\lfloor nt_{C(j)-1} \rfloor < i_1^j<i_{k_j}^j \leq \lfloor nt_{C(j)} \rfloor\\ j\in [ |C| ]}}  \prod_{j=1}^{|C|} \Bigg\{ c  \lfloor i_{1}^j-i_{k_{j-1}}^{j-1}\rfloor ^{\alpha -1}  \frac{ \nu^{k_j-1}}{\Gamma(\alpha(k_j-1 ))} \lfloor i_{k_j}^j-i_1^j\rfloor ^{\alpha(k_j-1) -1}\Bigg\}\\& \hspace{65pt}  \times \exp\big(-\lambda (i_{k_{|C|}}^{|C|} -i_1^1)/n\big)
    \\ &\sim n^{\alpha \sum_{i\in  [|C|] } k_i }       t_{C(|C|)}^{\alpha\sum_{i\in  [|C|] }k_i }   \iint \cdots \int_{Q_C^*}
  \prod_{i=1}^{|C|} 
 \Big\{ c ( x_{i}'-x_{i-1})^{\alpha-1} \frac{ \nu^{k_i-1} (x_{i}-x_{i}')^{\al(k_i-1)-1}}{\Gamma((k_i-1)\al)}\Big\} \\ & \hspace{65pt}\times
  \exp\big({-\lambda t_{C(|C|)} x_{|C|}}\big)
   \prod_{\ell\in [|C|]}dx_{\ell}'dx_{\ell}  \\ &=n^{\alpha \sum_{i\in  [|C|] } k_i }    t_{C(|C|)}^{\alpha\sum_{i\in  [|C|] }k_i }  A^*(C,\{k_i\}),
\end{align*}
where
\begin{align*}
 Q_C^*=\Big\{ x_i',x_i, x_j: x_0=0,  
\frac{ t_{C(i)-1}}{t_{C(|C|)}}
 <x_i'<x_i<\frac{t_{C(i)}}{t_{C(|C|)}}
  ,i\in  [|C|]  \Big\}. \end{align*}

Step 3) We will show that $X_t^*$ is continuous in probability.
This follows from the fact that when $C=\{2\}, k_1=2,$ 
\[ \mathrm{E}\big(  ( X_{t_{2}}^*- X_{t_{1}}^* )^2 \big)=2 t_{2}^{2\alpha }A^*(C,\{k_i\}) =2 t_{2}^{2\alpha }  \int_{t_{1}/t_{2}}^1 \int_{ t_{1}/ t_{2}}^{x_1} c(x_1')^{\al-1}\frac{\nu(x_1-x_1')^{\al-1}}{\Gamma(\al)}\exp(-\lambda t_{2} x_{1})  dx_1' dx_1 ,\]  and it converges to zero if $ t_2=t_1+h$ and $h\to0.$ 
It is trivial that $X_t^{(*,n)}$ has non-decreasing sample path. 
Therefore, by Theorem 3 in \cite{bingham_1971_limit},  $X_t^{(*,n)}$  converges weakly to $X_t^*$ in the space of $D[0,\infty)$ with Skorokhod ($J_1$) topology.
\\ \qed \\

{\bf {Proof of Corollary \ref{cor2}.}}\\
   The results are easily derived from the mgfs (\ref{mgf3}, \ref{mgf4}).  For $C\subset [m], k_i\in \mathbb{N}, i\in [|C|]$,  \[   \mathrm{E}( \prod_{i\in [|C|] }(X_{t_{C(i)}}^*-X_{t_{C(i)-1}}^*)^{k_i} ) =   t_{C(|C|)}^{\alpha\sum_{i\in  [|C|] }k_i } \Big(\prod_{i\in {[|C|]}} k_i! \Big)A^*(C,\{k_i\}).
\]
Especially, if $ C=\{ 1,2\}, k_1=k_2=1,$
\[ \mathrm{E}((X_{t_2}^*-X_{t_1}^*)X_{t_1}^*)= t_{2}^{2\alpha } \int_{t_1/t_2}^{1}\int_{0}^{t_1/t_2} (x_{2}-x_{1})^{\al-1}e^{-\lambda t_2x_2}  dx_1 dx_2.\]
\\ \qed \\

{{\bf {Proof of Proposition \ref{prop}.}}\\
  {\it i)} 
Let $\hat p(s,t), \hat p^*(s,t)$ be the Fourier transform of $p(x,t), p^*(x,t),$ respectively, i.e., \[ \hat p(s,t)= E(e^{{\rm i} s X_t}), \hspace{8pt} \hat p^*(s,t)= E(e^{{\rm i} s X^*_t}),\] then it is easy to see by Theorems \ref{theo1}, \ref{theo2}, and their proofs that

 \begin{align*} 
    \hat p(s,t)=m({\rm i}s,t)= 1+ 
  t a_T {\rm i}s + a_T  \sum_{k=1}^{\infty} ({\rm i}s)^{k+1} \nu^{k}    
 \int_{0}^t \int_{0}^x   \frac{e^{-\lambda y}y^{k\alpha-1}}{\Gamma(k\alpha)} dy    dx,\end{align*}
 
 \begin{align*}
\hat p^*(s,t)=m^*({\rm i}s,t)=1+\sum_{k=1}^{\infty}({\rm i}s)^k \nu^{k} \int_{0}^t \frac{e^{-\lambda y}y^{k\alpha-1}}{\Gamma(k\alpha)} dy  . \end{align*} 
  
  Let $\bar{\hat p}^*(s,r)$ be the Laplace transform of $\hat p^*(s,t),$ then we have, for s, r in a small neighborhood of zero,
  \begin{align*}
    \bar{\hat p}^*(s,r)&=\int_{0}^{\infty}  e^{-rt} \hat p^*(s,t) dt 
\\ & =\frac{1}{r}+\sum_{k=1}^{\infty} ({\rm i} s)^k \int_{0}^{\infty}  e^{-rt} \int_{0}^t \frac{e^{-\lambda y} y^{k\alpha -1}}{ \Gamma(k\alpha )}  dy  dt \\&=\frac{1}{r}+\sum_{k=1}^{\infty} \frac{({\rm i} s)^k}{r} (r+\lambda)^{-k\alpha}
      =\frac{1}{r(1- ({\rm i} s)(r+\lambda)^{-\alpha})},\end{align*}  
  from which we obtain
  \begin{equation} \label{eq1}
     ((r+\lambda)^{\alpha}-({\rm i} s) ) \bar{\hat p}^*(s,r)=\frac{(r+\lambda)^{\alpha}}{r} . \end{equation}
Since the Laplace transform of the tempered fractional derivative of a function is  (see e.g. \cite{  ALRAWASHDEH2017892, SABZIKAR201514})
\[ \mathcal{L} (    D_t^{\alpha,\lambda } g  )(r)  = \int_{0}^{\infty} e^{-rt} D_t^{\alpha,\lambda } g(t) dt=(r+\lambda)^{\alpha } \mathcal{L} g(r),\]
we invert the Laplace transform in (\ref{eq1}) to obtain
\[ D_t^{\alpha,\lambda }\hat p^*(s,t)   -({\rm i} s)  \hat p^*(s,t)=D_t^{\alpha,\lambda }(1)  ,\]
and invert the Fourier transform to obtain
\[    D_t^{\alpha,\lambda }p^*(x,t) +\frac{\partial p^*(x,t)}{\partial x}=\delta(x)   D_t^{\alpha,\lambda }(1). \] 
  Let $\bar {\hat p}(s,r)$ be Laplace transform of $\hat p(s,t),$ then we have, for $s,r$ in a small neighborhood of zero,
  \begin{align*}
    \bar {\hat p}(s,r)&=\int_{0}^{\infty}  e^{-rt} {\hat p}(s,t) dt 
\\ & =\frac{1}{r}+\frac{{\rm i} sa_T}{r^2}+a_T\sum_{k=1}^{\infty} ({\rm i}s)^k \int_{0}^{\infty} e^{-rt} \int_{0}^t  \int_0^{x} \frac{e^{-\lambda y} y^{k\alpha -1}}{ \Gamma(k\alpha )}  dy  dx dt \\&=\frac{1}{r}+\frac{a_T {\rm i}s}{r^2}+a_T\sum_{k=1}^{\infty} ({\rm i}s)^{k+1}\frac{(r+\lambda)^{-k\alpha}}{r^2}
      =\frac{1}{r}+ \frac{a_T {\rm i}s}{r^2}\bigg(\frac{1}{1-{\rm i}s(r+\lambda)^{-\alpha}}\bigg),\end{align*}  
  from which we obtain
  \begin{equation*}
     ((r+\lambda)^{\alpha}-{\rm i}s ) \bar {\hat p}(s,r)=\frac{(r+\lambda)^{\alpha}-{\rm i}s}{r}+\frac{a_T{\rm i}s}{r^2}(r+\lambda)^{\alpha} . \end{equation*}
 We invert the Fourier transform to obtain
 \[  D_t^{\alpha,\lambda }\hat p(s,t)- {\rm i}s \hat p(s,t)=  D_t^{\alpha,\lambda }(1)- {\rm i}s+ a_T {\rm i}s D_t^{\alpha,\lambda }( t) ,\]
and the Laplace transform  to obtain
\[    D_t^{\alpha,\lambda }p(x,t) +\frac{\partial p(x,t)}{\partial x}=\delta(x)   D_t^{\alpha,\lambda }(1)+  \delta'(x)-\delta'(x){a_T}   D_t^{\alpha,\lambda }( t) . \]  
\\ {\it ii)}  From (\ref{mgf3}), 
\begin{align*}
    m^*(s,t)&=e^{-\lambda t} \sum_{k=0}^{\infty}\frac{(st^{\alpha})^k}{\Gamma(k\alpha +1)}+\lambda \int_0^t e^{-\lambda y} \frac{(s y^{\alpha})^k}{\Gamma(k\alpha +1)}  dy\\&=e^{-\lambda t} E_{\alpha}( s t^{\alpha} )+\lambda \int_0^t e^{-\lambda y} E_{\alpha}(sy^{\alpha} ) dy,
\end{align*}
  where \[E_{\alpha}(z)= \sum_{k=0}^{\infty}\frac{z^k}{\Gamma(k\alpha +1)} \] is the Mittag-Leffler function.
It is known that the Laplace transform of inverse stable subordinator is the Mittag-Leffler function \citep{meerschaert_2019_stochastic}, therefore,
\begin{align*}
    m^*(s,t)=e^{-\lambda t} \int_{0}^{\infty} e^{sx} h(x,t) dx +\lambda \int_0^t e^{-\lambda y} \int_{0}^{\infty} e^{sx} h(x,y) dx    dy,
\end{align*}
and the result follows.
Also,
    \[ m(s,t)= a_T\int_0^t s m^*(s,y) dy+1,\] from which we can derive the following relationship between the marginal distribution of $X_t$ and $X_t^*.$
    \[ p(x,t)= - a_T \int_0^t \frac{\partial p^*(x,y)}{\partial x} dy +\delta(x). \]
 \\ \qed \\

{\bf {Proof of Proposition \ref{prop2}.}}\\
{\it i)} From Theorem \ref{theo1}, it follows that
\begin{align*}
    \mathrm{E}(\nabla X_t \nabla X_s)=&b  e^{-\lambda T} \int_{s-u}^s  \int_{t-u}^t e^{-\lambda (y-x)} dx dy  \\&=c'e^{-\lambda (s-t)}(e^{\lambda d}-2+e^{-\lambda d}  ),
\end{align*}
where $c'= b  e^{-\lambda T}/\lambda ^2.$\\
{\it ii)} Note that
\begin{align*}
\mathrm{E}(\nabla X_{s} \nabla X_{t})&=
\mathrm{E}( X_{ s-t}( X_{d}- X_{s-t} ))+\mathrm{E}( X_{ d}( X_{d+(s-t)}- X_{d} ))+\mathrm{E}( X_{ d-(s-t)}^2)
.\end{align*} By Theorem \ref{theo1} and the change of variable, 
\begin{align*}
    \mathrm{E}( X_{ s-t}( X_{d}- X_{s-t} ))&=
    b e^{-\lambda T}  \int_{s-t}^{d}\int_0^{s-t} e^{-\lambda (x_2-x_1)} dx_1 dx_2\\&= b e^{-\lambda T}/\lambda ^2(e^{-\lambda (s-t) }-e^{-\lambda d)})(e^{\lambda (s-t)}-1),
\end{align*}
and  similarly we obtain \[  \mathrm{E}( X_{d}( X_{d+(s-t)}- X_{d} )) = b e^{-bT}/\lambda ^2(e^{-\lambda d}-e^{-\lambda (d+ s-t)})(e^{\lambda d}-1).\] Also,  we have  \[E( X_{d-(s-t) }^2)=2b e^{-\lambda T}/\lambda ^2(e^{-\lambda (d-(s-t))}-1+\lambda (d-(s-t)
)).\]
Therefore
\begin{align*}
  &\mathrm{corr}(\nabla X_{t}, \nabla X_{s}) \\ &= \frac{b   e^{-\lambda T}/\lambda ^2 \big( e^{-\lambda (d+(s-t) )}+e^{-\lambda (d-(s-t)) }-2e^{-\lambda (s-t) } +2\lambda (d-(s-t) ) \big)  - (b e^{-\lambda T} d)^2 }{2b e^{-\lambda T}/\lambda ^2(e^{-\lambda d}-1+\lambda d)- (b e^{-\lambda T} d)^2  } , \end{align*}
and the result follows.
\\ \qed \\

{\bf {Proof of Proposition \ref{prop3}.}}\\
{\it i)} From Theorem \ref{theo1}, it is easy to see that when $s-t>d,$
\begin{align*}
\mathrm{E}(\nabla X_{t} \nabla X_{s})&= (s-t+d)^{\alpha+1}\frac{b  T^{\alpha-1}}{\Gamma(\alpha)^2}\int_{(s-t)/(s-t+d)}^1\int_{0}^{d/(s-t+d)} (x_2-x_1)^{\al-1} dx_1 dx_2\\& =\frac{b  T^{\alpha-1}}{\Gamma(\alpha+2)\Gamma(\alpha)} \big((s-t+d)^{\alpha+1}-2(s-t)^{\alpha+1}+ (s-t-d)^{\alpha+1}\big).
\end{align*}
{\it ii)} If $s-t<d,$
\begin{align*}
\mathrm{E}(\nabla X_{t} \nabla X_{s} )&=
\mathrm{E}( X_{s-t}( X_{d}- X_{s-t} ))+\mathrm{E}( X_{ d}( X_{d+(s-t)}- X_{d} ))+\mathrm{E}( X_{ d-(s-t)}^2),
\end{align*} and
\begin{align*}
    \mathrm{E}( X_{s-t}( X_{d}- X_{s-t} ))&=d^{\alpha+1}
    b T^{\alpha-1}/\Gamma(\alpha)^2  \int_{(s-t)/d}^{1}\int_0^{(s-t)/d} (x_2-x_1)^{\alpha-1} dx_1 dx_2\\&=\frac{b T^{\alpha-1}}{\Gamma(\alpha)\Gamma(\alpha+2)}\Big(d^{\alpha+1}-\big(d-(s-t)\big)^{\alpha+1}- (s-t)^{\alpha+1}\Big).
\end{align*}
Similarly, we obtain \[  \mathrm{E}( X_{ d}( X_{d+(s-t)}- X_{d} )) =\frac{b  T^{\alpha-1}}{\Gamma(\alpha)\Gamma(\alpha+2)}\Big(\big(d+(s-t)\big)^{\alpha+1}-d^{\alpha+1}-(s-t)^{\alpha+1}\Big).\] Also, we have  \[\mathrm{E}( X_{ d- (s-t)}^2)=\frac{2b T^{\alpha-1}  \big(d-(s-t) \big)^{\alpha+1}}{\Gamma(\alpha)\Gamma(\alpha+2)}.\]
Therefore, we have
\[  \mathrm{corr}( \nabla X_{t}, \nabla X_{s})  = \frac{c'\big( (d+(s-t))^{\alpha+1}-2(s-t)^{\alpha+1} +(d-(s-t))^{\alpha+1} \big)  - (b/\Gamma(\al) T^{\alpha-1} d)^2 }{2 c' d^{\alpha+1}- (b/\Gamma(\al) T^{\alpha-1} d)^2} , \] where
$c'=  b T^{\alpha-1}/ (\Gamma(\alpha)\Gamma(\alpha+2))$,
from which the result follows.
  \\ \qed\\

{\bf {Proof of Theorem \ref{theo3}.}}\\
We will prove the result for $Z_t$. For $Z_t^*$, it is proved similarly.

Let $t_0=0$. Note that
\begin{align*}
     \mathrm{E}(e^{{\rm i}s_1Z_{t_1}+ {\rm i}s_2 Z_{t_2}+\cdots + {\rm i}s_m Z_{t_m} })&=\mathrm{E}\Big(\exp\big( \sum_{i=1}^m {\rm i} S_i( Z_{t_i}- Z_{t_{i-1}}  ) \big) \Big)\\&=
    \mathrm{E}\Big[  \mathrm{E}\Big(\exp\big( \sum_{i=1}^m  {\rm i} S_i( Z_{t_i}- Z_{t_{i-1}}  ) \big)  \Big| \{ X_{t_i}, i\in [m] \}  \Big) \Big]\\
     &=\mathrm{E}\Big[  \prod_{i\in [m]}\mathrm{E}\Big(e^{{\rm i}S_iY_1} \Big)^{ X_{t_i} - X_{t_{i-1}}} \Big]
    \\&= \mathrm{E}\Big[    \prod_{i\in [m]}e^{\psi(S_i)(X_{t_i}- X_{t_{i-1}} )} \Big]\\&= m\Big( \big( \psi (S_1) ,\psi (S_2) , \cdots ,\psi (S_m)  \big), 
(t_1,t_2,\cdots, t_m)   \Big),
\end{align*}
where the last equality was derived from Theorem \ref{theo1} and its proof as follows. 
Since 
\begin{align*}
    \sum_{\substack{ k_i\in \{0\}\cup \mathbb{N}\\i\in[m]  }}   \big|{\psi'(S)}\big|^{\sum_{i\in [m]} k_i} 
\mathrm{E} \bigg( \prod_{i\in [m]} \frac{  (X_{t_i}- X_{t_{i-1}} )^{k_i} }{ k_i!} \bigg)   <\infty,
\end{align*}  where $\psi'(S)  =\max\{|\psi(S_i)|  : i\in [m]  \},$  therefore,
\begin{align*}
 \mathrm{E}\Big[    \prod_{i\in [m]}e^{\psi(S_i)(X_{t_i}- X_{t_{i-1}} )} \Big]&= \mathrm{E}\Big[    \prod_{i\in [m]} \sum_{k_i=0}^{\infty}{\psi(S_i)^{k_i}(X_{t_i}- X_{t_{i-1}} )^{k_i} /k_i!} \Big]\\&= \sum_{\substack{ k_i\in \{0\}\cup \mathbb{N}\\i\in[m]  }}   \Big( \prod_{i\in [m]}{\psi(S_i)}^{ k_i} \Big)
\mathrm{E} \bigg( \prod_{i\in [m]} \frac{  (X_{t_i}- X_{t_{i-1}} )^{k_i} }{ k_i!} \bigg)\\&=m\Big( \big( \psi (S_1) ,\psi (S_2) , \cdots ,\psi (S_m)  \big), 
(t_1,t_2,\cdots, t_m)   \Big)  . \end{align*} 
\\ \qed \\

{\bf {Proof of Corollary \ref{cor3}.}}\\
The results follow since
\begin{align*}
 \mathrm{E}(  \nabla Z_t \nabla Z_s   )& = \mathrm{E}\big( \mathrm{E}(  \nabla Z_t \nabla Z_s | \{  X_t, X_{t-d}, X_s, X_{s-d} \}  ) \big) \\&=\mathrm{E}\big(  \nabla X_t \nabla X_{s} \mathrm{E}(Y)^2  \big)=\mathrm{E}(Y)^2\mathrm{E}(\nabla X_t \nabla X_s) ,
\end{align*} for $t-s>d$, and in general,
\begin{align*}
 \mathrm{E}(  \nabla Z_{t_1} \nabla Z_{t_2}  \cdots \nabla Z_{t_k}  )& = \mathrm{E}\big( \mathrm{E}(  \nabla Z_{t_1} \nabla Z_{t_2} \cdots  \nabla Z_{t_k} | \{  X_{t_1}, X_{t_1-d},\cdots,  X_{t_k}, X_{t_k-d} \}  ) \big) \\&=\mathrm{E}\Big(  \prod_{i=1}^k \nabla X_{t_i} \mathrm{E}(Y)^k  \Big)=\mathrm{E}(Y)^k\mathrm{E}\Big(\prod_{i=1}^k \nabla X_{t_i} \Big) ,
\end{align*} 
for $t_i-t_{j}>d, i,j=1,2,\cdots, k. $
\\ \qed \\

\section*{Declaration}
The author did not receive support from any organization for the submitted work and  have no relevant financial or non-financial interests to disclose.

\bibliographystyle{abbrvnat}
\bibliography{main}

\end{document}